\documentclass{amsart}
\usepackage{graphicx}
\usepackage[title]{appendix}

\usepackage{amsmath}
\usepackage{amssymb}
\usepackage{amsthm}

\usepackage{algorithm}
\usepackage{algpseudocode}
\DeclareMathAlphabet{\mathcal}{OMS}{cmsy}{m}{n}
\usepackage{hyperref}
\usepackage{cleveref}
\crefname{claim}{claim}{claims}

\newtheorem{theorem}{Theorem}[section]
\newtheorem{proposition}{Proposition}[section]
\newtheorem{lemma}{Lemma}[section]

\newtheorem{claim}{Claim}

\theoremstyle{definition}
\newtheorem{definition}[theorem]{Definition}

\theoremstyle{remark}
\newtheorem{remark}[theorem]{Remark}

\numberwithin{equation}{section}
\newcommand{\FF}{\mathbb{F}}

\newcommand{\ttT}{\mathtt{T}}

\newcommand{\R}{\mathbb{R}}
\newcommand{\C}{\mathbb{C}}

\newcommand{\cR}{\mathcal{R}}
\newcommand{\cF}{\mathcal{F}}
\newcommand{\cB}{\mathcal{B}}

\newcommand{\tc}{\mathtt{c}}
\newcommand{\tb}{\mathtt{b}}

\newcommand{\cT}{\mathcal{T}}
\newcommand{\hcT}{\hat{\mathcal{T}}}
\newcommand{\hcB}{\hat{\mathcal{B}}}

\newcommand{\cG}{\mathcal{G}}
\newcommand{\cS}{\mathcal{S}}
\newcommand{\bLx}{\mathrm{L}_X}

\newcommand{\Lin}{\mathtt{Lin}}

\newcommand{\cM}{\mathcal{M}}

\newcommand{\frI}{\mathfrak{I}}

\newcommand{\cU}{\mathcal{U}}
\newcommand{\real}{\text{Re}}

\newcommand{\bA}{\mathrm{A}}

\newcommand{\bF}{\mathrm{F}}

\newcommand{\bG}{\mathrm{G}}
\newcommand{\bC}{\mathrm{C}}
\newcommand{\bCx}{\mathrm{C}_{\bx}}

\newcommand{\bI}{I}

\newcommand{\cN}{\mathcal{N}}

\newcommand{\ft}{\mathsf{T}}

\newcommand{\tK}{\mathtt{K}}

\newcommand{\bL}{\mathrm{L}}
\newcommand{\cL}{\mathcal{L}}

\newcommand{\fS}{\mathfrak{s}}

\newcommand{\bx}{X}

\newcommand{\bz}{z}
\newcommand{\bR}{\mathfrak{r}}

\newcommand{\blbd}{\lambda}
\newcommand{\blambda}{\vec{\lambda}}

\newcommand{\bP}{\mathrm{P}}
\newcommand{\hP}{\hat{\bP}}

\newcommand{\NCM}{\textsc{NCM}}
\newcommand{\ONCM}{\textsc{O-NCM}}

\newcommand{\EIn}{E}
\newcommand{\EE}{\mathcal{E}}

\newcommand{\EL}{\mathcal{E}_{L}}
\newcommand{\bLlbd}{\mathrm{L_{\blbd}}}
\newcommand{\bLlbdm}{\mathrm{L_{\blbd}^-}}

\DeclareMathOperator{\diag}{diag}
\DeclareMathOperator{\bdiag}{bdiag}

\DeclareMathOperator{\Imag}{Im}
\DeclareMathOperator{\Null}{Null}

\DeclareMathOperator{\grad}{grad}

\AtBeginDocument{
  \label{CorrectFirstPageLabel}
  
}
\begin{document}
  \title[Rayleigh Quotient cubic convergence]{Rayleigh Quotient Iteration, cubic convergence, and second covariant derivative}
\author{Du Nguyen}
\email{nguyendu@post.harvard.edu}
\address{organization={Indepdendent},
            city={Darien},
            postcode={06820}, 
            state={CT},
            country={USA}
}

\begin{abstract}
  {We generalize the Rayleigh Quotient Iteration (RQI) to the problem of solving a nonlinear equation where the variables are divided into two subsets, one satisfying additional equality constraints and the other could be considered as (generalized nonlinear Lagrange) multipliers. This framework covers several problems, including the (linear\slash nonlinear) eigenvalue problems, the constrained optimization problem, and the tensor eigenpair problem. Often, the RQI increment could be computed in two equivalent forms. The classical Rayleigh quotient algorithm uses the {\it Schur form}, while the projected Hessian method in constrained optimization uses the {\it Newton form}. We link the cubic convergence of these iterations with a {\it constrained Chebyshev term}, showing it is related to the geometric concept of {\it second covariant derivative}. Both the generalized Rayleigh quotient and the {\it Hessian of the retraction} used in the RQI appear in the Chebyshev term. We derive several cubic convergence results in application and construct new RQIs for matrix and tensor problems.}
\end{abstract}

\maketitle
\section{Introduction}
Several important problems in mathematics could be reduced to solving a vector equation of the form $\cL(\bx, \blbd) = 0$, where the function $\cL$ and the unknowns are divided into two groups as follows
\begin{equation}\label{eq:LMR0}
  \cL(\bx, \blbd) = \begin{bmatrix}\bL(\bx, \blbd) \\  \bC(\bx)\end{bmatrix}.\end{equation}
Here, $\bx$ and $\blbd$ are vector variables, defined on two vector spaces $\EE$ and $\EL$, respectively. The function $\bL(\bx, \blbd)=0$ in the first equation involves all variables, $\bL$ maps $\EIn\times \EL$ to $\EE$. In the second group, the {\it constraint} $\bC(\bx)=0$, involves only $\bx$, $\bC$ maps $\EE$ to the constraint vector space $\EL$. The vector variable $\blbd\in\EL$ plays the role of (generalized Lagrange) {\it multipliers}. Thus $\cL$ maps $\EIn\times \EL$ to itself. This setup covers at least four classes of problems encountered in the literature:

{\it a) The eigenvector problem:}
\begin{equation}
\begin{aligned}
\bL(\bx, \blbd) = & A \bx -\bx \blbd,\\
\bC(\bx) = & \frac{1}{2}(\bx^{\ft}\bx - 1)
\end{aligned}\label{eq:ceigen}
 \end{equation}
where $\bx$ is an $n\times 1$ vector, $A$ is an $n\times n$ matrix, $\lambda$ is a scalar. Here, $\EE$ is the space of $n\times 1$ vectors and $\EL$ is the base field, $\R$. A related problem has $\bL(\bx, \blbd) =A\bx - \bx\blbd - \tb= 0$ where $\tb\neq 0$ is a vector (\cite{Hazan} and references therein).

{\it b) The constrained optimization problem.}
This is the problem of optimizing a cost function $f$ in a space $E$ under the constraint $\bC(\bx)=0$ where $\bC$ is a function with target space $\EL$ and $\blbd\in \EL$ are the Lagrange multipliers. Let $\bC'$ be the Jacobian of $\bC$ then the Lagrangian multiplier equation is
$$\bL(\bx, \blbd) =\nabla f(\bx) - \bC'(\bx)^{\ft}\blbd=0.$$
Note $\bC'(\bx)^{\ft}$ is a map from $\EL$ to $\EE$. The system \cref{eq:LMR0} gives us the set of critical points.

{\it c) The nonlinear eigenvalue problem}:
\begin{equation}
\bL(\bx, \blbd) =\bP(\lambda)\bx = 0.
\end{equation}
In the real case, we set $\EL = \R$, with $\lambda\in \R$ is scalar. Here $\bP$ is a matrix with polynomial entries in $\lambda$. While this is not in the form of \cref{eq:LMR0} we can impose the constraint $\bC(\bx) = \bx^{\ft}\bx - 1$ (or $\bC(\bx) = \bz^{\ft}\bx-1$ for a fixed vector $\bz$). See \cite{guttel_tisseur_2017} for a survey.

{\it d) The tensor eigenpairs problem}: Here, $\EE=\FF^n$ is a vector space over a base field $\FF$, $\FF=\R$ or $\C$. Set $\EL=\FF$ , thus $\lambda$ is scalar. Let $\cT$ and $\cB$ be two vector-valued functions from $\EE$ to itself, with entries $\cT_i$ and $\cB_i, 1\leq i\leq n$ are homogeneous polynomials of degrees $m-1$ and $d-1$ respectively. The evaluations of $\cT$ and $\cB$ at $\bx\in \EE$ are written $\cT(\bx^{[m-1]})$ and $\cB(\bx^{[d-1]})$. Set 
\begin{equation}\label{eq:tensordef}
  \begin{aligned}
\bL(\bx, \blbd) = & \cT(\bx^{[m-1]})-\blbd\cB(\bx^{[d-1]}).
\end{aligned}
\end{equation}
A popular constraint is $\bC(\bx) = \frac{1}{2}(\bx^{\ft}\bx - 1)$ for the real case, and we will study $\bC(\bx) = \frac{1}{2}(\bx^*\bx - 1)$ for the complex case ($*$ is the Hermitian transpose). An important case is when $\cT=1/m\hcT'$, where $\hcT'$ denotes the gradient of a scalar homogeneous polynomial $\hcT$ of order $m$, and $\cB(\bx) = \bx$. These eigenpairs could be used to determine if $\hcT$ is nonnegative (\cite{LiqunQi}).

The Rayleigh quotient iteration (RQI) is among the most powerful methods to compute eigenvalues and vectors. For a vector $\bx$, the Rayleigh quotient is $\blbd = \frac{\bx^{\ft}A \bx}{\bx^{\ft}\bx}$, or $\bx^{\ft}A \bx$ on the unit sphere. In the $i$-th step, with $\lambda_i$ computed from $\bx_i$ by this equation, the iteration computes
\begin{equation} \bx_{i+1} = \frac{(A-\blbd_i\bI)^{-1}\bx_i}{||(A-\blbd_i\bI)^{-1}\bx_i||}.\label{eq:cray}\end{equation}
It has cubic convergence when $A$ is normal and quadratic otherwise (for suitable initial points).

Similar iterations have been suggested for the remaining problems. Points satisfying $\bC(\bx)=0$ are called feasible points. At the $i$-th step, where $\bx=\bx_i$ is feasible, a (vector-valued) function $\cR(\bx)$ of $\bx$ is used to compute $\lambda$ for the iterative step, then an intermediate step $\hat{\bx}_{i+1}$ is produced by solving a linear equation depending on $\lambda$, and a feasibility perturbation is applied to produce the next feasible step $\bx_{i+1}$. The feasibility perturbation is often available for common constraints, thus the overall procedure consists of constructing the generalized Rayleigh quotient $\cR$ and the iteration equation. In the literature, convergence results are obtained separately for each problem. We show here there is a common procedure and convergence analysis for these problems, including criteria for cubic convergence.

Recall an iteration $\{\bx_i\}$ converging to $\bx_*$ has order $k$ if $\lvert \bx_{i+1}-\bx_*\rvert \leq M\lvert \bx_{i}-\bx_*\rvert^k$ for some $M$, for large enough $i$. Quadratic and cubic convergence corresponds to $k=2$ and $k=3$. It is well-known Newton's method has quadratic convergence order and the Chebyshev method \cite{Nechepurenko} achieves cubic convergence. We link the cubic convergence of the Rayleigh quotient iteration with a {\it constrained Chebyshev term}. While very little geometric prerequisite is needed in the proof or in applications, we explain this Chebyshev term is essentially a {\it second covariant derivative} using {\it two connections} constructed from the problem formulation.

\subsection{Generalized RQI: quadratic convergence}\label{subsec:introRQI}
We denote the Jacobian of a nonlinear map $\bF$ between vector spaces $V$ and $W$ by $\bF'$. For $\bx\in V$, $\bF'(\bx)$ is a linear operator between $V$ and $W$. Let $\Lin(V, W)$ be the space of linear operators from $V$ to $W$. We will assume sufficient smoothness. Consider a constraint $\bC(\bx)=0$ defining a feasible set $\cM$ such that $\bC'(\bx)$ is of full rank in the subset of interest, $\cM$ is equipped with a feasibility perturbation $\bR$. Recall the tangent space of $\cM$ at $\bx\in\cM$ is the nullspace $T_{\bx}\cM$, or $T_{\bx}$ for short, of $\bC'(\bx)$. We are mainly interested in feasibility perturbation from the tangent space. For $\bx\in\cM, \eta\in T_{\bx}$, instead of considering a map from $\EE$ to $\cM$, mapping $\bx+\eta$ to a point on $\cM$, following \cite{AdlerShub}, we consider a map $\bR$ sending $(\bx, \eta)$ to $\bR(\bx, \eta)\in \cM$ for $\eta$ sufficiently small, with $\bR(\bx, 0) = \bx$, ($0=0_{T_{\bx}}\in T_{\bx}$ is the zero vector). For the constraint $\bx^{\ft}\bx =1$, we can use the rescaling $\bR(\bx, \eta) = \frac{1}{|\bx+\eta|}(\bx+\eta)=\frac{1}{(1+|\eta|^2)^{1/2}}(\bx + \eta).$

Fix $\bx$, then $\bR(\bx, .): \eta\mapsto\bR(\bx, \eta)$ is a map from $T_{\bx}$ to $\EE$. Assume it has the Taylor expansion
    \begin{equation}\bR(\bx, \eta) = \bx +\eta + O(|\eta|^2)\label{eq:retract}
    \end{equation}
    or $\bR$ is a {\it retraction}, thus, the partial derivative $\bR_{\eta}$ of $\bR$ in $\eta$, considered as a linear map from $T_{\bx}$ to $\EE$ satisfies $\bR_{\eta}(\bx, 0)\eta=\eta$. This holds for the rescaling perturbation of the sphere and most feasibility perturbations in practice, and we can normalize a reasonable perturbation to this form.

The second ingredient of a generalized RQI is a projection function. In \cref{eq:ceigen,eq:cray} of the classical RQI, with $\lambda = \bx^{\ft}A\bx$
    $$\bF(\bx):=\bL(\bx, \bx^{\ft}A\bx) = A\bx - \bx\bx^{\ft}A\bx = (I-\bx\bx^{\ft})A\bx.$$
    The matrix $\Pi(\bx)= I-\bx\bx^{\ft}$ is a {\it projection}, as $(I-\bx\bx^{\ft})^2=I-\bx\bx^{\ft}$. Thus $\Pi(\bx)\bF(\bx)=\Pi(\bx)^2A\bx= \bF(\bx)$. This seems incidental, but it is crucial in generalizing the RQI.
  
    An operator $P\in\Lin(\EE, \EE)$ is called an affine projection if $P^2=P$. It is often constructed in two ways. If $A^-$ is a left inverse of a linear operator $A$, $A^-A = I$, then $AA^-$ is a projection.  If $P$ is a projection then $I_{\EE} - P$ is also a projection. If $\eta\in\EE$ is in the range $\Imag(P)$ of $P$ then $P\eta = \eta$. In the special case where $P$ is self-adjoint, $P$ is called an {\it orthogonal projection}. We call a smooth map $\Pi$ from $\cM$ to $\Lin(\EE, \EE)$, the space of linear operators on $\EE$ an (affine) {\it projection function}, or just {\it projection} when there is no confusion, if, for $\bx\in\cM$, $\Pi(\bx)$ is an affine projection.  

{\it We will assume for a projection function $\Pi$, $\Pi(\bx)$ is of rank $\dim\EE-\dim\EL$ for $\bx\in\cM$.}

The following claim is the key to our generalized RQI to solve \cref{eq:LMR0}.
\begin{claim}\label[claim]{claim:cl1}Consider a smooth feasible set $\cM$ defined by $\bC(\bx) = 0$, $\bC'(\bx)$ is of full rank for $\bx\in \cM$, and $\Pi:\cM\to\Lin(\EE, \EE)$ is a projection function. Let $\bF$ be a smooth function from $\cM$ to $\EE$ satisfying
\begin{equation}\Pi(\bx)\bF(\bx) =\bF(\bx).\end{equation}
For an initial point $\bx_0\in \cM$ sufficiently close to a solution $\bx_*$ of $\bF(\bx)=0$, the iteration $\bx_{i+1} = \bR(\bx_i, \eta_i)$ with $\bR$ satisfies \cref{eq:retract} and $\eta_i\in T_{\bx_i}\cM=T_{\bx_i}$ solving the system (\ref{eq:lmr10}), (\ref{eq:lmr2}) below
\begin{gather}
  \Pi(X_i)\bF'(X_i)\eta_i = -\bF(X_i)\label{eq:lmr10}\\
  \bC'(X_i)\eta_i = 0\label{eq:lmr2}
\end{gather}
is well-defined and converges quadratically to $\bx_*$ if the system is non-degenerated at $\bx_*$.
\end{claim}
To use this claim to solve \cref{eq:LMR0}, note for any projection $\Pi$, $\bL(\bx, \lambda) = 0$ is equivalent to
\begin{gather}\Pi(\bx)\bL(\bx, \lambda) = 0,\label{eq:zeroOfSection}\\
  \bL(\bx, \lambda)-  \Pi(\bx)\bL(\bx, \lambda) = 0.\label{eq:section}
\end{gather}
If a function $\cR(\bx)$ substituted to $\lambda$ satisfying \cref{eq:section} for all $\bx$, consider $\bF(\bx) := \bL(\bx, \cR(\bx))$, then $\Pi(\bx)\bF(\bx) =\bF(\bx)$ is satisfied by design, and we can use the iteration in the claim.

To simplify \cref{eq:lmr10}, we construct $\Pi$ as below. Let $\bLlbd(\bx, \blbd)$ be the partial derivative of $\bL$ in $\lambda$, considered as a linear map from $\EL$ to $\EE$. Assume it is onto, which is the case when the inverse function theorem is satisfied for \cref{eq:LMR0}, then it has a left inverse $\bL_{\lambda}^-(\bx, \blbd)$ (any computationally efficient inverse or the Moore-Penrose inverse can be used). Solve for $\lambda=\cR(\bx)$ from
\begin{equation}\bL_{\lambda}^-(\bx, \blbd)\bL(\bx, \blbd) = 0.\label{eq:Rayleigh}\end{equation}
Then with $\Pi^- = \bLlbd\bL_{\lambda}^-$, $\Pi(\bx):= I_{\EE} - \Pi^-(\bx, \cR(\bx))$ is a projection, \cref{eq:section} becomes \cref{eq:proj1} and \cref{eq:proj2} holds as $\bLlbd =\bLlbd\bL_{\lambda}^-\bLlbd$ by the left inverse assumption
  \begin{gather}\Pi(\bx)\bL(\bx, \cR(\bx)) =\bL(\bx, \cR(\bx)),
\label{eq:proj1}\\
    \Pi(\bx)\bLlbd(\bx, \cR(\bx))=0\label{eq:proj2}.
\end{gather}    
From \cref{eq:proj2}, the expression $\Pi(\bx)\bF'(\bx)$ in \cref{claim:cl1} with $\bF(\bx)=\bL(\bx, \cR(\bx))$ reduces to
  $$\Pi(\bx)\bF'(\bx)=\Pi(\bx)\{\bLx(\bx, \cR(\bx)) + \bLlbd(\bx, \cR(\bx))\cR'(\bx) \} = \Pi(\bx)\bLx(\bx, \cR(\bx)).$$
Here $\bLx(\bx, \lambda)$ is the partial derivative in $\bx$. We will choose $\cR$ solving \cref{eq:Rayleigh}.

For the eigenvalue problem $\bL(\bx, \lambda) = A\bx - \bx\lambda$ with constraint $\bx^{\ft}\bx = 1$, $\bLlbd(\bx, \lambda) = -\bx$, take the left inverse $\bL_{\lambda}^- = -\bx^{\ft}$, then \cref{eq:Rayleigh} gives us $\lambda = \bx^{\ft} A\bx$ and $\Pi(\bx) = \bI_{\EE} - \bx\bx^{\ft}$. Using the iteration in \cref{claim:cl1} for $\bF(\bx) = A\bx -\bx\bx^{\ft}A\bx$, we get the Rayleigh quotient iteration, since
\begin{lemma}[Schur form]\label{lem:cl2}Let $P\in\Lin(\EE, \EE)$ be an affine projection, $B\in\Lin(\EE, \EE)$ be an invertible operator, $F\in\EE$, and $D\in\Lin(\EE, \EL)$. Consider the linear equations in $\eta$
  \begin{equation}\label{eq:schurEq}
    \begin{gathered}P B\eta = -P F,\\
      D\eta =0.
      \end{gathered}
  \end{equation}
If $H\in\Lin(\EL, \EE)$ is such that $P H = 0$ and $DB^{-1}H$ is invertible, $\eta$ below solves \cref{eq:schurEq}
\begin{equation}\label{eq:schurEq1}
  \eta = -B^{-1}F +B^{-1}H(DB^{-1}H)^{-1}DB^{-1} F.
  \end{equation}
\end{lemma}
The proof is by direct substitution, see \cref{appx:SchurEq}.

Continuing with the eigenvalue problem, with $\lambda = \bx^{\ft}A\bx$, $\bF(\bx) = A\bx - \bx\bx^{\ft}A\bx$, $\bLx(\bx,\blbd)\eta = (A - \lambda I_{\EE})\eta$ for $\eta\in\EE$. In the lemma, set $B = A - \lambda I_{\EE}$,  $H:\lambda \mapsto \bx\lambda$, $D:\eta\mapsto \bx^{\ft}\eta$, $F= A\bx - \bx\lambda$, then $B^{-1}F = \bx$ and $\zeta:=B^{-1}H =(A-\lambda I)^{-1}\bx$, hence
$$\eta = -\bx + \zeta(\bx^{\ft}\zeta)^{-1}\bx^{\ft}\bx= -\bx + (\bx^{\ft}\zeta)^{-1}(A-\lambda I)^{-1}X$$
with $(\bx^{\ft}\zeta)^{-1}$ is a scalar. Then $\bx +\eta$ is proportional to $(A-\lambda I)^{-1}X$, and the retraction gives us the iteration  $X_{i+1} =\|(A-\lambda I)^{-1}X_i\|^{-1}(A-\lambda I)^{-1}X_i$.

This ``explains''  why the Rayleigh quotient iteration has quadratic convergence and relates it to the Newton-type iteration of \cref{claim:cl1}. This explanation for the classical Rayleigh quotient essentially appeared in \cite[example 5]{AdlerShub}. Note $\Pi$ is hidden in the classical RQI.

In \cref{lem:cl2}, take $H=\bLlbd(\bx, \cR(\bx))$, we have the Schur form for \cref{claim:cl1}.  We will apply this result to the remaining three problems in \cref{sec:examples}. The Schur form's usage is restricted by the condition $\bLx(\bx, \cR(\bx))$ is invertible. In optimization problems (\cref{subsec:opt}), \cref{eq:lmr10} is often solved directly. Let $Q_{\Pi}$ be an orthonormal basis of the nullspace of $\bL_{\lambda}^-(\bx, \blbd)$ (same as the range of $\Pi(\bx)$), and $Q_T$ an orthonormal basis of $T_{\bx}$, then expressing $\bL$ and $\Pi\bLx$ in these bases (evaluated at $(\bx, \blbd)$) we have
\begin{equation}\eta = -Q_T(Q_{\Pi}^{\ft} \Pi \bLx Q_T)^{-1}Q_{\Pi}^{\ft}\bL\label{eq:NewtonForm},\end{equation}
the {\it Newton form} of the iteration. Note $Q_{\Pi}Q_{\Pi}^{\ft}$ and $\Pi(\bx)$ are both projections to the range of $\Pi(\bx)$, hence, $Q_{\Pi}Q_{\Pi}^{\ft}\Pi = \Pi, Q_{\Pi}Q_{\Pi}^{\ft}\bL=\bL$ which helps verify $\Pi\bLx\eta = Q_{\Pi}Q_{\Pi}^{\ft}\Pi\bLx\eta = -Q_{\Pi}Q_{\Pi}^{\ft}\bL=-\bL$.

To summarize, we have a procedure to solve \cref{eq:LMR0} by a generalized RQI

1. Find a left inverse $\bL_{\lambda}^{-}(\bx, \blbd)$ of $\bL_{\lambda}(\bx, \blbd)$ such that $\bL_{\lambda}^{-}(\bx, \blbd)\bL(\bx, \blbd)=0$ is easy to solve for $\lambda$ as a function $\lambda = \cR(\bx)$ in $\bx$. Choose a retraction $\bR$.

2. Determine if the Schur form or the Newton form solution of \cref{eq:lmr10} is preferable. For the Newton form, also set up $\Pi(\bx)=I_E - \bL_{\lambda}(\bx, \blbd)\bL_{\lambda}^{-}(\bx, \lambda)$ with $\lambda=\cR(\bx)$.

3. Apply \cref{alg:alRQI}.

\begin{algorithm}
\label{alg:RQI}
\begin{algorithmic}
\State{Input: $\bx_0\in\cM$, choice of retraction $\bR$, Rayleigh quotient $\cR$ satisfying \cref{eq:proj1,eq:proj2}}.\;
\For{$i=0,1,\cdots$}\;
\State Compute $\blbd_{i} = \cR(\bx_i)$\;
\State Compute $\eta_i\gets -\bLx^{-1}\bL + \bLx^{-1}\bLlbd(\bC'\bLx^{-1}\bLlbd)^{-1}\bC'\bLx^{-1}\bL$ at $\bx_i, \lambda_i$;\Comment{Schur form}
\State \quad or solve $\Pi(\bx_i)\bLx(\bx_i, \lambda_i)\eta_i = -\bL(\bx_i, \lambda_i)$ using \cref{eq:NewtonForm};\Comment{Newton form}
\State Compute $\bx_{i+1} \gets \bR(\bx_i, \eta_i)$\;\Comment{Terminal condition is verified after this step}
\EndFor
\end{algorithmic}
\caption{Generalized Rayleigh quotient iteration}
\label{alg:alRQI}
\end{algorithm}
\subsection{Geometric interpretation}
\label{sec:geometry}
The discussion in this section and \cref{sec:geomcubic} are not required to follow the rest of the paper. However, we believe they are interesting interpretations of the results in geometric terms, clarifying several geometric concepts that are often presented more abstractly.

The projection function  $\Pi$ on the feasible set $\cM$ associates to each feasible point a vector space $\EE_{\bx} :=\Imag(\Pi(\bx))$. The collections $\EE_{\Pi}=\cup_{\bx\in \cM}\EE_{\bx}$ and $T\cM=\cup_{\bx\in \cM}T_{\bx}$ are  {\it vector bundles} in the geometric literature, assuming appropriate smoothness and constant rank of $\EE_{\bx}$ and $T_{\bx}$. They could be considered smooth subsets of $\EE^2$ of pairs $(\bx, \omega)$ with $\bx$ satisfying a nonlinear constraint $\bC(\bx) =0$ not involving $\omega$, while
the constraint on $\omega$ is linear. We have $\EE_{\Pi} = T\cM$ in the eigenvalue problem, but they are different for the RQI in \cref{rem:Btensor}.

A function $\bF$ from $\cM$ to $\EE$ satisfying the condition $\Pi(\bx)\bF(\bx)=\bF(\bx)$ is called a {\it section} of the bundle $\EE_{\Pi}$, we require $\bF(\bx)\in \EE_{\bx}$.
Our RQI framework could be understood as a construction of a projection function $\Pi$ defining a vector bundle $\EE_{\Pi}$, (and the {\it connection $\nabla$} below) and a section $\bF(\bx) = \bL(\bx, \cR(\bx))$. A section $\tc$ of the tangent bundle (a function to $\EE$ with $\tc(\bx)\in T_{\bx}$) is called a {\it vector field}, each $\tc(\bx), \bx\in\cM$ is a tangent vector at $\bx$.
  
For a section $\bF$ and a tangent vector $\eta\in T_{\bx}, \bx\in\cM$, we define
\begin{equation}\nabla_{\eta}\bF(\bx):=\Pi(\bx)\bF'(\bx;\eta)  = \bF'(\bx; \eta) - \Pi'(\bx, \eta))\bF(\bx).
\end{equation}
Here, $\Pi'(\bx, \eta)$ and $\bF'(\bx, \eta)$ are directional derivatives of $\Pi$ and $\bF$ in direction $\eta$, the second equality follows by differentiating $\Pi(\bx)\bF(\bx) = \bF(\bx)$ in direction $\eta$. If $\tc$ is a vector field then $\nabla_{\tc}\bF:\bx\mapsto \Pi(\bx)\bF'(\bx;\tc(\bx))$ is a section of $\EE_{\Pi}$. The expression for $\nabla_{\tc}\bF$ shows it is a {\it covariant derivative}, or {\it connection} in the differential geometric sense \cite[section 5.2]{AMS_book}. Intuitively, a covariant derivative is a rule to take directional derivative of sections, resulting in sections. The Newton increment equation reads $\nabla_\eta\bF(\bx) = -\bF(\bx)$.

Thus, \cref{claim:cl1} is a {\it Newton method on a vector bundle}, see \cite{AdlerShub,Wang,FERREIRA2002} and references therein. Newton method on the tangent bundle is well-studied in the Riemannian optimization literature. The novelty here is the identification of a generalized Rayleigh quotient with the construction of a vector bundle and a connection from the data of \cref{eq:LMR0} using a left-inverse, even for nonlinear multipliers, and the convergence analysis with a retraction (see \cite{AdlerShub,AMS_book} for the tangent bundle case), versus geodesics in \cite{Wang}. The Schur form iteration for the general case is also new.

\subsection{Cubic convergence and Chebyshev iterations}
For normal matrices, the classical RQI has cubic convergence and the two-sided iterations in \cite{Ostrowski1959,schreiber_thesis} also converge cubically for any matrix. The second goal of the paper is to clarify the condition for cubic convergence. We will define notations more properly in \cref{subsec:notation}, but briefly, the Hessian of a function $\bF$ could be considered as a map valued in bilinear functions, we use the notation $\bF^{(2)}(\bx;\eta^{[2]}) = \bF^{(2)}(\bx)\eta^{[2]}$ to denote the evaluation of this map at $\bx$, and of the bilinear function in both linear variables at $\eta$. The (partial) Hessian $\bR_{\eta\eta}(\bx, 0; \eta^{[2]})$ of the retraction $\bR$ in the variable $\eta$ at $(\bx, 0)$ and the directional derivative $\cR'(\bx; \eta)$ also appear in the analysis.
\begin{claim}\label[claim]{claim:cl3}For a function $\bF$ satisfying $\Pi(\bx)\bF(\bx) = \bF(\bx)$ as in \cref{claim:cl1}, and for a tangent vector $\eta\in T_{\bx}$ at $\bx$, set
  $$\begin{gathered}
  \bG(\bx; \eta^{[2]}) = \bF^{(2)}(\bx; \eta^{[2]}) + \bF'(\bx)\bR_{\eta\eta}(\bx, 0; \eta^{[2]}).
  \end{gathered}$$
  If $\Pi(\bx_*)\bG(\bx_*; \eta^{[2]}) = 0$ for all $\eta\in T_{\bx_*}$ at a solution $\bx_*$, the RQI in \cref{claim:cl1} converges cubically to $\bx_*$ when $\bx_0$ is sufficiently close to $\bx_*$. If $\eta_i = \bx_i^{\cN_T}$ is the RQI increment in \cref{claim:cl1}, then the Rayleigh-Chebyshev iteration $\bx_{i+1} = \bR(\bx_i, \bx_i^{\cN_C})$ with $\bx_i^{\cN_C}=\eta_i -\frac{1}{2}\bx_i^{\tau} \in T_{\bx}$ and $\bx_i^{\tau}\in T_{\bx_i}$ satisfies
  $$  \Pi(\bx_i)\bF'(\bx_i)\bx_i^{\tau} = \Pi(\bx_i)\bG(\bx_i)\eta_i^{[2]}$$
  converges cubically for an initial point $\bx_0$ sufficiently close to $\bx_*$.\hfill\break
If $\bF = \bL(\bx, \cR(\bx))$, $\lambda=\cR(\bx)$ arises from an RQI as in \cref{subsec:introRQI}, we can use $\bG_{\bL}$ in place of $\bG$:
  $$\begin{gathered} \bG_{\bL}(\bx; \eta^{[2]}) = \bL_{\bx\bx}(\bx, \lambda)\eta^{[2]} + 2\bL_{\bx\blbd}(\bx, \lambda)[\eta, \cR'(\bx; \eta)]+
    \bL_{\blbd\blbd}(\bx, \blbd)\cR'(\bx; \eta)^{[2]}  + \bLx(\bx, \lambda)\bR_{\eta\eta}(\bx, 0; \eta^{[2]}).
  \end{gathered}
  $$
\end{claim}
The above is a constrained version of the Chebyshev iteration, see \cref{sec:HOGaussNewton}. For the eigenvalue problems, expanding $\bR(\bx, \eta) =\frac{\bx+\eta}{(1+|\eta|^2)^{1/2}}$ to Taylor series in $\eta$ to the second term gives $\bR_{\eta\eta}(\bx, 0; \eta^{[2]}) =-\bx|\eta|^2$. In \cref{sec:examples}, we use this to verify cubic convergence for several eigenvalue RQIs. The term $\bF'(\bx)\bR_{\eta\eta}(\bx, 0; \eta^{[2]})$ often vanishes at a solution $\bx_*$ for homogeneous tensor\slash eigenvalue problems. In \cref{sec:geomcubic}, we relate it to a {\it connection} constructed from $\bR$, and give an interpretation of the Chebyshev term as a {\it second covariant derivative}.

We apply the general analysis here to construct new RQIs for the tensor eigenpair problem in \cref{subsec:tensor}. A unitary version of the Schur form RQI over $\C$ finds all complex eigenpairs and also identifies the real pairs as a by-product in the generic case, improving on previous work \cite{JaffeWeissNadler}. We find new complex pairs for the Motzkin polynomial  \cite[example 5.9]{CARTWRIGHT2013942}, completing the eigenpair count for this classical example. We verify cubic convergence for several eigenvalue problems, propose a cubic convergence iteration for the generalized eigenvalue problem, and derive Chebyshev iterations in a few interesting cases. A natural question is if we can develop a framework for constrained homotopy continuation \cite{Homotopy}, using RQI-type iterations. We hope this work will lead to further research along this line.

\subsection{Notations and outline}\label{subsec:notation}
We mimic the convention of \cite{OrtegaRhein} for derivatives. We work with a base field $\FF$ which is $\R$ or $\C$. We denote by $\FF^{n\times m}$ the space of $n\times m$ matrices on $\FF$, by $\Lin(V, W)$ the space of linear map between vector spaces $V$ and $W$. The zero vector in $V$ is denoted $0_V$, or just $0$. The inner product of two (column) vectors $\bx_1, \bx_2$ is $\bx_1^{\ft}\bx_2$ or $\bx_1^{*}\bx_2$, identifying $\FF^{1\times 1}$ with $\FF$.

For a vector-valued function $\bF$ from an open subset $D\subset V$ to $W$, by $\bF'$ we denote the Jacobian or Fr{\'e}chet derivative, and the directional derivative at $\bx\in D$ in direction $\eta$ is written $\bF'(\bx)\eta$ or $\bF'(\bx; \eta)$, the latter is preferable when we need grouping.  Thus, if $D$ is an open subset in $V$, $\bF'$ is an operator-valued function, $\bF'(\bx)\in \Lin(V, W)$. Higher (partial) derivatives could be considered as a map from $D$ to the space of multilinear maps from $V$ to $W$. Thus, at $\bx\in D$, the $l$-th-order derivative of $\bF$, denoted by $\bF^{(l)}(\bx)$ is a $l$-linear map from $V$ to $W$. For $\eta_1,\cdots,\eta_l\in V$, we denote its evaluation as $\bF^{(l)}(\bx; \eta_1,\cdots \eta_l)$ or $\bF^{(l)}(\bx; [\eta_1,\cdots \eta_l])$. In the expressions $\bF^{(l)}(\bx) \eta^{[l]}=\bF^{(l)}(\bx; \eta^{[l]})$, $\eta^{[l]}$ denotes $\eta$ being repeated $l$ times. If $\Phi$ is an operator value function from $V$ to $\Lin(W_1, W_2)$, with $V, W_1, W_2$ are vector spaces and $\Phi(\bx)\omega$ denotes its valuation at $\bx\in V$ operating on $\omega\in W_1$, then $\Phi^{(l)}(\bx; \eta_1,\cdots \eta_l)\omega$ denotes the $l$-th order derivative of $\Phi$ evaluated in direction $\eta_1,\cdots \eta_l$ and operate at $\omega$. In general, we use the round brackets for base variable evaluations or for groupings of base and directional variables. The semicolumn separates the base variable(s) from the directional derivative variables. The $l$-terms multi-dimensional Taylor series expansion around $X_0$ is \cite[NR 3.3-3]{OrtegaRhein}
\begin{equation}\begin{gathered}
    \bF(\bx) = \bF(\bx_0) + \sum_{j=1}^l \frac{1}{j!}\bF^{(j)}(\bx_0;[\bx - \bx_0]^{[j]}) + \int_0^1  \frac{(1-t)^l}{l!}\bF^{(l+1)}((1-t)\bx_0 + t\bx; [\bx-\bx_0]^{[l+1]})dt    
\end{gathered}\end{equation}

Partial derivatives are sometimes denoted by the usual subscript convention, eg $\bLx, \bLlbd$ when convenient. However, we will use a position-based notation for more complex formulas. If $\bF$ is defined using two vector variables $\bx, Y\in D_1\times D_2$ for two sets $D_1, D_2$, then we write $\bF(\bx, .)$ for the map $Y\mapsto \bF(\bx, Y)$ with fixed $\bx$ in the appropriate domain, and write $\bF(., Y)$ similarly. The partial derivatives in $Y$ evaluated at $Y$ are $\bF(\bx, .)'(Y)=\bF_Y(\bx, Y)$, which is a linear operator, and $\bF(\bx, .)'(Y; \eta)$ denotes its valuation at $\eta$. Thus, for a function $\bF$ and a retraction $\bR$, the directional derivative in direction $(0, \xi)$ of $\bF(\bR(\bx, \eta))$  is $\bF'(\bR(\bx, \eta); \bR(\bx, .)'(\eta; \xi))$ by the chain rule. The notation $\bR(\bx, .)'(\eta)\xi=\bR(\bx, .)'(\eta;\xi)$ is also used.

Denote by $A_{V_0\downarrow W_0}$ the restriction of a map $A\in \Lin(V, W)$ from two spaces $V$ and $W$ to subspaces $V_0\subset V, W_0\subset W$, with $A(V_0)\subset W_0$. The inverse, if exists, is denoted by $A_{V_0\downarrow W_0}^{-1}$, shorthand for $(A_{V_0\downarrow W_0})^{-1}$.

{\it To focus on the main ideas, we will assume sufficient smoothness, and use the word smooth for short}. Most results hold for class $C^3$ or $C^4$. 

The main idea of the convergence analysis is to require the retraction $\bR$ maps an open ball $B_{\bx_f}(\rho)$ of the tangent space of a feasible point $\bx_f$ (an initial or final point) to an open subset of $\cM$ in \cref{sec:retract}, then translate the problem to one of linear constraint, with a nontrivial retraction $\fS$. The linear constraint case in \cref{sec:HOGaussNewton} is the most technical, the proofs are in the appendices. The iterations in \cref{sec:RQIRC} are as described in the introduction, with applications in \cref{sec:examples}. The relationship with the second covariant derivative is in \cref{sec:geomcubic}.

\section{Feasibility perturbation and retractions}\label{sec:retract}
We state the definition and prove the basic properties of a retraction in this section, illustrated with the example of the rescaling retraction of the sphere. We also define the rescaling retractions and compute their Taylor expansion in the tangent variable.
\subsection{Definition and basic properties}
Assume $\bC:\EE\to \EL$ is a smooth map, with Jacobian $\bC'(\bx)$ which is {\it onto} $\EL$ at any $\bx\in \EE$ with $\bC(\bx) = 0$. Let $\cM$ be the corresponding feasible set, the solution set of $\bC(\bx) = 0$. This is a multidimensional smooth surface, a {\it manifold}. Denote by  $T\cM$, its tangent bundle, the subset of $\EE^2 = \EE\times \EE$ of pairs $(\bx, \eta)$ such that
\begin{equation}\label{eq:eqTTM}
  \begin{gathered}\bC(\bx) = 0,\\
    \bC'(\bx)\eta=0.
    \end{gathered}
\end{equation}
For a fixed $\bx\in\cM$, the vector space defined by $\bCx(\bx)\eta=0$ is the tangent space at $\bx$, denoted by $T_{\bx}\cM$, or $T_{\bx}$ for short subsequently, and $\eta$ is called a tangent vector. In elementary physics, an element $\eta\in T_{\bx}\cM$ is considered as a velocity vector of a particle moving smoothly on $\cM$. The Jacobian of the constraints \cref{eq:eqTTM} is a block diagonal matrix with two diagonal blocks $\bC'(\bx)$, hence is {\it onto} $\EL^2$, thus $T\cM$ is also a smooth manifold. In the next two paragraphs, we summarize the basic geometric facts required, well-known for surfaces and curves.

Locally, $\cM$ is parametrized by open subsets of $E_T := \R^{\dim\EE-\dim\EL}$. For $\bx\in\cM$, a {\it local parametrization} around $\bx$ is defined as a  smooth map $\phi$, injective, from an open subset $\Omega$ of $E_T$ into $\EE$, with $\bx\in\phi(\Omega)\subset \cM$, such that $\phi^{-1}$, defined on $\phi(\Omega)$, is continuous, and $\phi'(z)$ is injective for $z\in \Omega$. The full rank assumption of $\bC'(\bx)$ implies there is a local parametrization around $\bx$.

Differentiating $\bC(\phi(z)) = 0$ in direction $v\in \EE_T$, we get $\bC'(\phi(z))\phi'(z)v = 0$, or $\phi'(z)v$ is tangent to $\cM$ at $\phi(z)$. The dimension count and the injective assumption of $\phi'(z)$ imply $\phi'(z)$ is bijective to $T_{\phi(z)}$. Let $\Omega_2$ be an open subset of $\EE_T$, then for $(z, v)\in \Omega\times\Omega_2\subset \Omega\times \EE_T$, define
\begin{equation}\label{eq:Dphi}
  D_{\phi}:(z, v)\mapsto (\phi(z), \phi'(z)v)\in T\cM
\end{equation}
then $D_{\phi}$ is a local parametrization of $T\cM\subset \EE^2$.

A map $\bR$ from an open subset $T_1$ of $T\cM$ to $\cM$ is smooth at $(\bx, \eta)\in T_1$ if for a local parametrization $D_{\phi}$ on $\Omega\times\Omega_2$ around $(\bx, \eta)$ and $\phi_1$ on $\Omega_1$ around $\bR(\bx, \eta)$, the combined map $\phi^{-1}_1\circ\bR\circ D_{\phi}$ from $D_{\phi}^{-1}(T_1)\cap\Omega\times\Omega_2$ to $\Omega_1$ is smooth. This is independent of parametrizations.

We are interested in a perturbation $\bR$ mapping a pair $(\bx, \eta)\in T\cM$ to $\cM$ for $\eta$ sufficiently small, with $\bR(\bx, 0) = \bx$, since our Newton increments are tangent vectors.
\begin{lemma}\label{lem:vectorTrans}
  Fix $\bx\in\cM$. If $\tK$ is a $C^1$ map from a neighborhood $\cU\subset T_{\bx}$ of $\xi\in T_{\bx}$ to $\cM\subset\EE$ then for $\eta\in T_{\bx}$, $\tK'(\xi)$ as a map from $T_{\bx}$ to $\EE$ satisfies $\tK'(\xi)\eta\in T_{\tK(\xi)}$, or $\tK'(\xi)(\cU)\subset T_{\tK(\xi)}$.
\end{lemma}
\begin{proof} Differentiate $\bC(\tK(\xi+t\eta)) = 0$ in $t$ at $t=0$, we get $\bC'(\tK(\xi))\tK'(\xi; \eta) = 0$ by the chain rule.
\end{proof}

Applying this to feasibility perturbations, for $\xi\in T_{\bx}$, then $\bR(\bx, .)'(\xi)$ could be considered as a map from $T_{\bx}$ to $T_{\bR(\bx, \xi)}$. In particular, $\bR(\bx, .)'(0)$ maps $T_{\bx}$ to itself. Note $\bR(\bx, .)(\xi)'\in\Lin(T_{\bx}, \EE)$, but we will denote $\bR(\bx, .)'(\xi)^{-1}$ the inverse, if exists, of $\bR(\bx, .)'(\xi)_{T_{\bx}\downarrow T_{\bR(\bx, \xi)}}$ (recall this means restricting the range of $\bR(\bx, .)'(\xi)$ to the image).

\begin{remark}
  Technically, it is simplest to consider a retraction as a map from $T\cM$ onto $\cM$ satisfying certain requirements on the Taylor series. This works for the rescaling map on the sphere, as $(X+\eta)/\lvert X+\eta\rvert$ is defined for any $\eta\in T_X$. However, we do want to consider, for example, rescaling maps for other constraints where rescaling only works for tangent vectors close enough to $0_X$. The following makes explicit the requirement that the radius where a perturbation exists is not too small, required in iterations, otherwise, it is adapted from \cite{AdlerShub,AbM}.  
\end{remark}

\begin{definition}\label[definition]{def:BR}  For a point $\bx_f\in\cM$, a {\it retraction} around $\bx_f$ (or a retraction to $\Omega^{\bR}$) is a smooth map $\bR$ from an open subset $\ttT_{\bR}\subset T\cM$ to $\cM$, with range in $\cM$ containing a neighborhood $\Omega^{\bR}$ of $\bx_f$, and there is $\rho_{\bR} > 0$ such that $\ttT_{\bR}$ contains
  \begin{equation}\label{eq:TOmegaR0}\ttT\Omega^{\bR}_{\rho_{\bR}} := \{(\bx, \eta)|\; (\bx, \eta)\in T\cM, \bx\in \Omega^{\bR}, |\eta|< \rho_{\bR}\}.\end{equation}
If $\bR(\bx, 0)$ is defined then we require $\bR(\bx, 0) = \bx$ and $\bR(\bx, .)'(0) = I_{T_{\bx}}$.
\end{definition}
\begin{remark}\label{rem:BRlocal}Let $\phi$  be a local parametrization of $\cM$ with $\phi(0) = \bx_f$. Requiring the existence of $\Omega^{\bR}$ and $\rho_{\bR}$ such that $\bR$ is defined in $\ttT\Omega^{\bR}_{\rho_{\bR}}$ is equivalent to requiring the range of $\bR$ containing $\{(\phi(z)|\; |z| \leq \rho_1\}$ and the domain of $\bR$ contains a set $\ttT\Omega^{\bR}_{\phi, \rho_1, \rho_2}$ for some $\rho_1, \rho_2> 0$ where
\begin{equation}\ttT\Omega^{\bR}_{\phi, \rho_1, \rho_2} := \{(\phi(z), \phi'(z)v)|\; |z| < \rho_1, |v| < \rho_2\}\label{eq:TOmegaR}
\end{equation}
This follows as $\phi'(z)$ and $\phi'(z)^{-1}$ exist and are continuous for a parametrization, thus bounded uniformly on a bounded set, so we can translate a bound on $|\eta| =|\phi'(z)v|$ to a bound on $v$, and vice versa.

We can normalize a feasibility perturbation to a retraction. If $\bR$ is feasibility perturbation satisfying the conditions of the retraction except for $\bR(\bx, .)'(0)$ is only assumed to be invertible, then $\hat{\bR}(\bx, \eta) := \bR(\bx, (\bR(\bx, .)'(0))^{-1}\eta)$ is a retraction since $\hat{\bR}(\bx, .)'(0; \eta) = \bR(\bx, .)'(0; (\bR(\bx, .)'(0))^{-1}\eta)=\eta.$
\end{remark}

If $\cM$ is the unit sphere, for a feasible point $\bx_f$, the rescaling of $\bx_f+\eta$ for $\eta\in T_{\bx_f}$ maps $T_{\bx_f}$ bijectively to the hemisphere centered at $\bx_f$, and this map parametrizes the hemisphere.

We now generalize this result, showing $\bR(\bx_f, .)$ parametrizes an open set around $\bx_f$.
\begin{proposition}\label[proposition]{prop:retractionTopo}Assume $\bR$ is a retraction around a feasible point $\bx_f\in\cM$.
  
  1) There is a radius $\rho > 0$ such that the map $\tK_f:=\bR(\bx_f, .):\eta\mapsto\bR(\bx_f, \eta)$ from the ball $B_{\bx_f}(\rho) :=\{\eta\in T_{\bx_f}||\eta| < \rho\}$ to its image $C_{\bx_f}(\rho)\subset\cM$ is a local parametrization. The ball $B_{\bx_f}(\rho)\subset T_{\bx_f}$ is called a {\it retraction ball}, $\tK_f(B_{\bx_f}(\rho)) = C_{\bx_f}(\rho)\subset\cM$ is called the {\it retraction cap}.

  2) Thus, the set $\ttT\Omega^{\bR}_{\tK_f, \rho_1, \rho_2}$ in \cref{eq:TOmegaR} exists, for each $\bx_f\in \Omega$, there are radii $\rho_1 >0 $, $\rho_2 > 0$ such that $\tK_f$ maps $B_{\bx_f}(\rho_1)$ to $C_{\bx_f}(\rho_1)$ with invertible Jacobian, and for all $\bx\in C_{\bx_f}(\rho_1)$, and $\delta\in T_{\bx_f}$ with $|\delta| < \rho_2$, if $\bx = \bR(\bx_f, \xi)$ then $\bR(\bx, \tK_f'(\xi; \delta))$ and $\tK_f^{-1}\bR(\bx, \tK_f'(\xi; \delta))$ exist.
\end{proposition}
\begin{proof}Let $\phi$ be a local parametrization of $\cM$ around $\bx_f$. By \cref{rem:BRlocal}, there are radii $r_0, \Delta_0$ such that $\phi(B_0)$ is in the range of $\bR$, where $B_0$ is the ball radius $r_0$ in $\EE_T = \R^{\dim\EE - \dim\EL}$, and $\ttT\Omega^{\bR}_{\phi, r_0, \Delta_0}$  exists.

  Since $\tK_f := \bR(\bx_f, .)$ is continuous, $\tK_f^{-1}(\phi(B_0))$ is open in $T_{\bx_f}$, thus there is a neighborhood $\cU\subset  T_{\bx_f}$ of $0_{T_{\bx_f}}$ such that $\tK_f(\cU)\subset \phi(B_0)$. Consider the map $h: \cU\times \EE_T\to\EE_T$, $h(\eta, z) = \phi^{-1}(\tK_f(\eta)) - z$ for $(\eta, z)\in  \cU\times \EE_T$. Since $\tK_f'(0)\xi=\xi$, by the chain rule, the partial derivative $h_{\eta}(0, 0)$ is $\xi\mapsto \phi'(0)^{-1}\xi\in T_{\bx_f}$ for $\xi\in T_{\bx_f}$, where $ \phi'(0)$ is invertible as a map to $T_{\bx_f}$. Thus, the implicit function theorem \cite[5.2.4]{OrtegaRhein} applies, we have open balls $B_1\subset B_0$ centered at $0_{\EE_T}$, $B_{\eta}\subset T_{\bx_f}$ at $0_{T_{\bx_f}}$, and a unique function $\eta_f(z)$ with $h(\eta_f(z), z) = 0$ for $z \in B_0, \eta_f(z)\in B_{\eta}$.

We have $\phi^{-1}(\tK_f(\eta_f(z)) = z$ or $\tK_f(\eta_f(z))=\phi(z)$. If $Y\in\phi(B_1)$, then $\eta = \eta_f\circ \phi^{-1}(Y)$ satisfies $\tK_f(\eta) = \phi (\phi^{-1}(Y)) =Y$. Since $\tK_f$ is continuous, $\tK_f^{-1}(\phi(B_1))$ is open, and its intersection with $B_1$ is open, containing a ball of radius $\rho$. Restricting to $B_{\bx_f}(\rho)$, the uniqueness of the implicit function shows $\tK_f$ is one-to-one, with continuous inverse $\eta_f\circ \phi^{-1}$. The remaining properties of parametrizations are verified. This proves 1).

Item 2) follows from \cref{rem:BRlocal}.
\end{proof}

We offer the first clue why the (partial) Hessian of $\bR$ appears in \cref{claim:cl3}.
\begin{proposition}\label[proposition]{prop:HessianIntrinsic}
  If $A$ is a smooth function on $\EE$ with values in a vector space $V$, assume $\bx\in \cM$ and $\eta\in T_{\bx}\cM$ then $A'(\bx)\eta$ is only dependent on the values of $A$ on $\cM$. If $\bR$ is a retraction then
  $$A^{(2)}(\bx; \eta^{[2]}) + A'(\bx)\bR(\bx, .)^{(2)}(0; \eta^{[2]})$$
is also only dependent on the values of $A$ on $\cM$.
\end{proposition}
\begin{proof}The first statement is well-known. Let $\phi$ be a local parametrization around $\bx$ defined on $\Omega\subset \EE_T$, and $\eta = \phi'(0)v$ for some $v\in \EE_T$. Consider the function $g(t) =A(\phi(tv))$ for $t$ sufficiently small. Then $\dot{g}(0) = A'(\phi(0))\phi'(0)v= A'(\bx)\eta$, but $g$ depends on values of $A$ in $\phi(\Omega)\subset\cM$ only.

If $\bR$ is a retraction, consider $f(t) = A(\bR(\bx, t\eta))$. The second derivative at $0$ of $f$ is
$$\ddot{f}(0) =\frac{d}{dt}_{|t=0}(A'(\bR(\bx, t\eta))\bR(\bx_0, .)'(t\eta; \eta)) = A^{(2)}(\bx; \eta^{[2]}) + A'(\bR(\bx, t\eta))_{t=0}\bR(\bx, .)^{(2)}(0, \eta^{[2]})
$$
by the chain rule. On the other hand, $A(\bR(\bx, t\eta))$ is dependent on the values of $A$ on $\cM$ only.
\end{proof}

We can get similar statements for higher derivatives. The function $\fS(\xi, \delta)$ below expresses a retraction to $\bx = \bR(\bx_f, \xi)\in \cM, \xi\in T_{\bx_f}$ in terms of a retraction to a chosen point $\bx_f\in \cM$. It appears when we translate an iteration on $\cM$ to an iteration on $T_{\bx_f}$ by a change of variable.
\begin{proposition}\label[proposition]{prop:fSretract}For $\bx_f\in\subset\cM$, set $\tK_f:= \bR(\bx_f, .)$, then with $\rho_1, \rho_2 >0$ as in 2) \cref{prop:retractionTopo},
\begin{equation}    
\fS(\xi, \delta; \bx_f)=  \fS(\xi, \delta) := \tK_f^{-1}\bR(\tK_f(\xi); \tK_f'(\xi; \delta))\in T_{\bx_f} \text{ for }\xi, \delta \in T_{\bx_f}\label{eq:FS}
\end{equation}
exists for $|\xi| < \rho_1$ and $|\delta| < \rho_2$ and is a retraction to $T_{\bx_f}\in \EE$. Moreover,
\begin{equation}
  \fS(\xi, \delta) = \xi+\delta +\frac{1}{2}(\tK_f'(\xi))^{-1}[
    \bR(\tK_f(\xi), .)^{(2)}(0; [
      \tK_f'(\xi; \delta)]^{[2]})  - \tK_f^{(2)}(\xi; \delta^{[2]})] + O(|\delta|^3)\label{eq:fSbR}
\end{equation}
\end{proposition}
\begin{proof}Recall for $\phi_1, \phi_2\in T_{\bx}$, $\bR(\bx, .)'(\phi_1; \phi_2)$ is $\lim_{t\to 0}\frac{1}{t}(\bR(\bx, \phi_1 +t\phi_2) -\bR(\bx, \phi_1))$. Rewrite the equation for $\fS$ then differentiate it with respect to $\delta$ in direction $\epsilon\in T_{\bx_f}$ twice
\begin{gather}  \tK_f(\fS(\xi, \delta)) = \bR(\tK_f(\xi); \tK_f'(\xi; \delta))\label{eq:calcFS0}\\
  \tK_f'(\fS(\xi, \delta); \fS(\xi, .)'(\delta; \epsilon)) = \bR(\tK_f(\xi), . )'(\tK_f'(\xi; \delta); \tK_f'(\xi; \epsilon))\label{eq:calcFS1}\\
        \tK_f^{(2)}(\fS(\xi, \delta); [\fS(\xi, .)'(\delta; \epsilon)]^{[2]}) +
\tK_f'(\fS(\xi, \delta); \fS(\xi, .)^{(2)}(\delta; \epsilon^{[2]}))
    = \bR(\tK_f(\xi), . )^{(2)}(\tK_f'(\xi; \delta); [\tK_f'(\xi; \epsilon)]^{[2]}).\label{eq:calcFS2}
\end{gather}
    Set $\delta$ to zero, we get the first three terms of the power series expansion of $\fS(\xi, \epsilon)$,
$$\begin{gathered}
\fS(\xi, 0) =  \tK_f^{-1}\bR(\bR(\bx_f, \xi), 0) =\bR(\bx_f, .)^{-1}\bR(\bx_f, \xi) = \xi,\\
\fS(\xi, .)'(0)[\epsilon] = (\tK_f'(\xi))^{-1}\bR(\tK_f, \xi), . )'(0; \tK_f'(\xi; \epsilon))=\epsilon,\\
\tK_f'(\xi; \fS(\xi, .)^{(2)}(0; \epsilon^{[2]}))  = \bR(\tK_f, \xi), . )^{(2)}(0; [\tK_f'(\xi; \epsilon)]^{[2]}) -    \tK_f^{(2)}(\xi; \epsilon^{[2]}).
\end{gathered}    
$$
    where we used $\bR(\tK_f, \xi), . )'(0; \tK_f'(\xi; \epsilon)) =\tK_f'(\xi; \epsilon)$ as $\bR$ is a retraction, note $(\tK_f'(\xi))^{-1}$ exists in the retraction ball. The statement follows from the Taylor series expansion.
\end{proof}
\subsection{The rescaling retraction on the unit sphere}
The computation of $\tK_f'$ for the unit sphere is known in \cite[section 8.1]{AMS_book}.
\begin{proposition}Consider the sphere $\cM=S^{n-1}\in \R^n$ ($n\geq 2$),  defined by the equation $X^{\ft}X = 1$. For $\bx_f \in S^{n-1}$, the retraction $\tK_f(\xi) = \bR(\bx_f, \xi) = \frac{1}{(1+ |\xi|^2)^{1/2}}(\bx_f +\xi)$ defined on $T_{\bR} = TS^{n-1}$ maps $\xi\in T_{\bx_f}$ to the hemisphere $\{\bx\in S^{n-1}|\;\bx_f^{\ft}\bx > 0\}$. For $\rho > 0$, let $C_{\bx_f}(\rho)$ be the retraction cap, then
  \begin{gather}
\label{eq:retSphereScale}
\frac{1}{(1+|\xi|^2)^{1/2}}(\bx+\xi)= \bx + \xi  - \frac{1}{2}\|\xi\|^2X + O(\|\xi\|^3),\\
C_{\bx_f}(\rho) = \{\tK_f(\xi) | |\xi| < \rho, \xi\in T_{\bx_f}\}
= \{\bx\in S^{n-1}| 1\geq \bx_f^{\ft}\bx > \frac{1}{(1+\rho^2)^{1/2}}\} \label{eq:Csphere},\\
    \tK_f^{-1}(\bx) =\bR(\bx_f, .)^{-1}\bx =\frac{1}{\bx_f^T\bx}\bx - \bx_f\quad\text{ for }\bx\in C_{\bx_f}(\rho), \rho > 0
\label{eq:Csphere0},\\
\tK_f'(\xi)\delta = \frac{1}{(1+|\xi|^2)^{1/2}}(\delta -\frac{\xi^{\ft}\delta}{1+|\xi|^2}(\bx_f+\xi))\text{ for }\delta\in T_{\bx_f},\\
   \fS(\bx_f; \xi, \delta) =
 \xi + \frac{1+|\xi|^2}{1+|\xi|^2 - \xi^{\ft}\delta}\delta=
 \xi + \delta + \frac{\xi^{\ft}\delta}{1 + |\xi|^2-\xi^{\ft}\delta}\delta\label{eq:fSsphere}.
  \end{gather}
  Any $\rho_1, \rho_2> 0$ satisfying condition \ref{eq:condrho12} below satisfies the first condition of 2), \cref{prop:retractionTopo}, that is $\fS(\bx_f; \xi, \delta)$ exists if $|\xi| < \rho_1, |\delta| < \rho_2$.
  \begin{equation}\label{eq:condrho12}
    \rho_2 < 2 \text{ or } (\rho_1 < 1 \text{ and } 1+|\rho_1|^2 -\rho_1\rho_2 > 0).
\end{equation}
\end{proposition}
\begin{proof} The Taylor series expansion of $(1+|\xi|^2)^{-1/2}$ gives \cref{eq:retSphereScale}. Since $\bx_f^{\ft}\bR(\bx_f, \xi) = \frac{1}{(1+|\xi|^2)^{1/2}} > 0$, and the function $(1+\rho^2)^{-1/2}$ is decreasing in $\rho$, the image of $\tK_f$ is in the hemisphere $\bx_f^T\bx > 0$ and (\ref{eq:Csphere}) follows. We can verify directly $\tK_f^{-1}\bx$ on the right-hand side of \cref{eq:Csphere0} satisfying $\bx_f^{\ft}(\tK_f^{-1}\bx) = 0$, and $\bx_f + \tK_f^{-1}\bx =\frac{1}{\bx_f^{\ft}\bx}\bx$ is proportional to $\bx$, thus $\frac{1}{|\bx_f + \tK_f^{-1}\bx|}(\bx_f + \tK_f^{-1}\bx) =\bx$.

A routine derivative gives $\tK_f'$. Since $\bR(\tK_f(\xi), \tK_f'(\xi, \delta))$ is proportional to $\tK_f(\xi)+ \tK_f'(\xi, \delta)$, we assume the scaling factor is of the form $c(1+|\xi|^2)^{1/2}$, then using $\bx_f^{\ft}\xi = \bx_f^{\ft}\delta = 0$
  $$\begin{gathered}
  \bR(\tK_f(\xi), \tK_f'(\xi, \delta)) =
c(1+|\xi|^2)^{1/2}(\tK_f(\xi)+ \tK_f'(\xi, \delta))=
  c(\bx_f +\xi +
    \delta -\frac{\xi^{\ft}\delta}{1+|\xi|^2}(\bx_f+\xi)),\\
    \bx_f^{\ft}\bR(\tK_f(\xi), \tK_f'(\xi, \delta)) = c(1 - \frac{\xi^{\ft}\delta}{1+|\xi|^2}),\\
    \tK_f^{-1}(\bR(\tK_f(\xi), \tK_f'(\xi, \delta))) =(1 - \frac{\xi^{\ft}\delta}{1+|\xi|^2})^{-1}(\bx_f +\xi +
    \delta -\frac{\xi^{\ft}\delta}{1+|\xi|^2}(\bx_f+\xi)) -\bx_f\\
    =(\bx_f+\xi) +(1 - \frac{\xi^{\ft}\delta}{1+|\xi|^2})^{-1}\delta -\bx_f = \xi + \frac{1+|\xi|^2}{1+|\xi|^2 - \xi^{\ft}\delta}\delta
  \end{gathered}$$
which gives us the formula for $\fS$. It exists if $f_1 = 1+|\xi|^2 -\xi^{\ft}\delta > 0$. If $|\delta|=\rho$ is constant, then $f_1$ is smallest if $\delta = \frac{\rho}{|\xi|}\xi$, with $f_1 = 1+|\xi|^2 - \rho|\xi|$. The region defined by $1 +\rho_1^2 -\rho_1\rho_2 > 0, \rho_1 > 0, \rho_2 > 0$ could be divided into three subregions, the union of two is described in \cref{eq:condrho12}, characterized by the condition that if $0 < z_1 < \rho_1, 0 < z_2 < \rho_2$ then $1 +z_1^2 -z_1z_2 > 0$. The remaining subregion with $\rho_1 \geq 1, \rho_2 \geq 2$ does not have this property to ensure $\fS(\xi, \delta)$ exists, if $|\xi| = z_1, |\delta| = \rho=z_2$.
\end{proof}

\subsection{Rescaling retractions}\label{subsec:rescaleProject}
We generalize the rescaling retraction on the sphere to feasible sets with a single constraint.
\begin{proposition}[Rescaling retraction]Assume $\EL = \R$, and $\bC'(\bx; \bx) \neq 0\in\EL=\R$ near $\bx_f\in\cM$. There is a neighborhood $\Omega^{\bR}\subset\cM$ of $\bx_f$ and $\rho_{\bR} > 0$ such that the equation in $\gamma\in \R$
\begin{equation}\bC(\gamma(\bx+\eta))=0\end{equation}
 for $(\bx, \eta)\in \ttT\Omega^{\bR}_{\rho_{\bR}}\subset T\cM$ (in \cref{eq:TOmegaR0}) has an implicit function solution with $\gamma=1$ at $\eta=0$. Then $\bR(\bx, \eta) = \gamma(\bx+\eta)$ is a retraction defined on $\ttT\Omega^{\bR}$. We have the Taylor expansion
  \begin{equation}\label{eq:scaleRetract}
    \bR(\bx, \eta) = X + \eta  -\frac{\bC^{(2)}(\bx; \eta^{[2]})}{2\bC'(\bx; \bx)}X + O(|\eta|^3).
  \end{equation}
\end{proposition}
\begin{proof}Let $\phi$ be a local parametrization from $\Omega_T\subset\EE_T=\R^{\dim(\EE)-\dim(\EL)}$ of $\cM$ near $\bx_f$. Consider $h(z, v, \gamma) = \bC(\gamma(\phi(z) + \phi'(z)v))\in\EL$ defined on $(z, v, \gamma)\in \Omega_{T}\times E_T\times \R$. Then $h(0, 0, 1) = \bC(\bx_f)=0$, $h_{\gamma}(0, 0, 1) = \bC'(\bx_f; \bx_f)\neq 0$. Thus, the implicit function theorem applies, giving us a ball $D=D_1\times D_2$ in $\Omega_T\times \EE_T$ of points $\{(z, v)\}$ with $|z| < \rho_1, |v| < \rho_2$ and $\delta_{\gamma} > 0$ such that $\gamma$ can be solved uniquely as a function $\gamma_{\phi}(z, v)$ with $\gamma_{\phi}(0, 0) = 1$ for $(z, v)\in D$ and $|\gamma - 1| < \delta_{\gamma}$. The conditions of \cref{rem:BRlocal} for $\bR(\bx, \eta) = \gamma(\bx+\eta)$ follows by properties of the implicit function. By local uniqueness of the implicit solution, $\gamma_{\phi}(z, 0) = 1$, or $\bR(\bx, 0) =\bx$ for $\bx\in \phi(D_1)$.

  Composing with $\phi^{-1}$, $\gamma=\gamma_{\phi}\circ \phi^{-1}$ is a function of $\bx$ and $\eta$, its partial derivatives are evaluated by implicit function rule for $g(\gamma, \bx, \eta) = \bC(\gamma(\bx+\eta))$. Hence, we can compute $\bR(\bx, .)'$ using
  $$\begin{gathered}\gamma(\bx, .)'(\eta; \xi) = -(\frac{\partial g(\gamma, \bx, \eta)}{\partial\gamma})^{-1}\frac{\partial g(\gamma, \bx, \eta)}{\partial\eta} = -\frac{\bC'(\gamma(\bx+\eta); \gamma\xi)}{\bC'(\gamma(\bx+\eta); \bx+\eta)},\\
  \bR_{\eta}(\bx, 0)\xi = \bR(\bx, .)'(0; \xi)= \gamma(\bx, .)'(0; \xi)(\bx+0) + \gamma(\bx, 0)\xi  = - \frac{\bC'(\bx; \xi)}{\bC'(\bx; \bx)}\bx+\xi=\xi\end{gathered}$$
  where $\xi\in T_{\bx}$, hence $\bC'(\bx, \xi) = 0$. Consider a fixed $\bx$, set $\gamma_{\eta}= \gamma(\bx, .)', \gamma_{\eta\eta}=\gamma(\bx, .)^{(2)}$. Take directional derivative of $\bC(\gamma(\bx+\eta))=0$ in $\eta$ in tangent directions $\xi_1, \xi_2$ consecutively
  $$\begin{gathered}\bC'(\gamma(\bx+\eta); \gamma_{\eta}(\eta; \xi_1)(\bx+\eta) + \gamma\xi_1 )=0,\\
  \bC^{(2)}(\gamma(\bx+\eta); \gamma_{\eta}(\eta; \xi_2)(\bx+\eta) + \gamma\xi_2, \gamma_{\eta}(\eta; \xi_1)(\bx+\eta) + \gamma\xi_1 )\\
  + \bC'(\gamma(\bx+\eta); \gamma_{\eta\eta}(\eta; \xi_2,\xi_1)(\bx+\eta) +
\gamma_{\eta}(\eta; \xi_1)\xi_2 + \gamma_{\eta}(\eta; \xi_2)\xi_1 )
    =0.\end{gathered}$$
    At $\eta = 0, \gamma=1$, $\gamma_{\eta}(0; \xi) = 0$ for $\xi\in T_{\bx}\cM$ since $\bC'(\bx, \xi)=0$ as before, the above gives us 
    $$\begin{gathered}\gamma_{\eta\eta}(0; \xi_1,\xi_2) = -\frac{\bC^{(2)}(\bx; \xi_1, \xi_2)}{\bC'(\bx; \bx)}.
    \end{gathered}$$
Thus, equation (\ref{eq:scaleRetract}) follows from $\bR(\bx, .)^{(2)}(0; \eta^{[2]}) =\gamma_{\eta\eta}(0; \eta^{[2]})\bx$ and the Taylor formula.
\end{proof}

The requirement $\bC'(\bx; \bx) \neq 0$ means the vector $\bx$ is not tangent to $\cM$. Note that when $\bC(\bx) = f(\bx) -1$ and $f$ is a homogeneous function of degree $n \neq 0$, Euler's identity implies $\bC'(\bx; \bx) = n f(\bx) = n\neq 0$, the equation $\bC(\gamma(\bx+\eta))=0$ has a simple solution $\gamma = f(\bx+\eta)^{-1/n}$ if $\eta$ is sufficiently small. The sphere corresponds to $f(\bx) = \bx^{\ft}\bx$ and $n=2$. In general, $\gamma$ is a solution of a scalar equation.

\section{Linear constraints and unconstrained case with a projection}\label{sec:HOGaussNewton}
This section is the most technical in the paper. The main proofs are deferred to the appendix.

Using the result of \cref{sec:retract}, we translate the constrained problem to one of linear constraints, by considering a reference point $\bx_f$, and solve $\cF(\xi) = \bF(\bR(\bx_f, \xi)) = 0$ instead of $\bF(\bx) = 0$, for $\xi\in T_{\bx_f}$ with the projection $\Phi(\xi) = \Pi(\bR(\bx_f,\xi))$,  $\Phi(\xi)\cF(\xi) =\cF(\xi)$, and the retraction $\fS$ in \cref{prop:fSretract}. The feasible set is now called $T_0$ ($=T_{\bx_f}$). The variable becomes $\xi$ instead of $\bx$.

A local linear constrained problem becomes an unconstrained problem if we parametrized the constraint set. Hence, we can consider the unconstrained problem over a vector space $T_0$.

Thus, we will start with two vector spaces $T_0$ and $\EE$,  with a function $\cF$ from an open subset $\Omega$ from $T_0$ to $\EE$, with Jacobian $\cF'$. We assume $\Omega$ to be convex, and $\dim T_0\leq \dim \EE$.

When $\dim T_0=\dim\EE$, for $\xi_*\in\Omega$, if $\cF(\xi_*) = 0$ with $\cF'(\xi_*)$ invertible, an inverse function $\cF_{\frI}$ exists, $\cF_{\frI}(\cF(\xi)) = \xi$ near $\xi_*$. For $Y\in \cU\subset\EE$ where the inverse function is defined, the $k$ terms Taylor expansion of $\cF_{\frI}$ at $\cF(\xi_i)$ (assumed to be in $\cU$) is
$$\cF_{\frI}(Y) = \xi_i + \sum_{j=1}^{k-1}\frac{1}{j!}\cF_{\frI}^{(j)}(\cF(\xi_i); [Y - \cF(\xi_i)]^{[j]}) + \int_0^1\frac{(1-t)^{k-1}}{(k-1)!}\cF_{\frI}^{(k)}(\cF(\xi_i)+ t(Y-\cF(\xi_i)); [Y - \cF(\xi_i)]^{[k]})dt.
$$
The Taylor coefficients for $\cF_{\frI}$ could be computed by the implicit-valued function theorem, or by power series substitution. We have $\cF_{\frI}'(Y)\eta = (\cF'(\xi))^{-1}\eta$ with $\xi = \cF_{\frI}(Y)$ and
$$\cF_{\frI}^{(2)}(Y; \eta^{[2]}) = -(\cF'(\xi))^{-1}\cF^{[2]}(\xi; [(\cF'(\xi))^{-1}\eta, (\cF'(\xi))^{-1}\eta ]).$$
Following \cite{Nechepurenko}, we can construct an iteration with order-$k$ convergence as follows. For $Y = 0$, $\cF_{\frI}(0)=\xi_*$, and if we set $\xi_{i+1}= \xi_i + \sum_{j=1}^{k-1}\frac{1}{j!}\cF_{\frI}^{(j)}(\cF(\xi_i); [ - \cF(\xi_i)]^{[j]})$,
the residual $\xi_* - \xi_{i+1}$, from the above expansion is $O(|\cF(\xi_i)|^k)= O(|\xi_i-\xi_*|^k)$ if we consider a bounded region and $\cF$ has bounded derivatives. Thus, $k$-th order convergence follows from the Taylor series of $\cF_{\frI}$ by design. The case $k=2$ corresponds to the Newton method with increment $\delta_2 = -(\cF'(\xi))^{-1}\cF(\xi)$, and the case $k=3$ corresponds to the Chebyshev method, with step
$\xi^{\bC} = \xi + \delta_2 -\frac{1}{2}(\cF'(\xi))^{-1}\cF^{(2)}(\xi; \delta_2^{[2]})
$.

When $\dim T_0< \dim\EE$ and $\cF'(\xi_*)$ is not invertible, assume there is an affine projection $\Phi$ such that $\Phi(\xi)\cF(\xi) = \cF(\xi)$ for $\xi\in\Omega$. We also assume $\dim(\Imag(\Phi(\xi)))=\dim(T_0)$. In \cref{prop:prelimOverdeter}, we define a left inverse $\cF^{\natural}$ of $\cF'$ and compute derivatives of $\cF_{\diamond}$, playing the role of $\cF_{\mathfrak{I}}$ (compare \cref{eq:bJFdiamZ,eq:jFDia2} with $\cF_{\frI}'$ and $\cF_{\frI}^{(2)}$). The proof is in \cref{appx:proofOverdeter}.

Recall the notation $A_{V_0\downarrow W_0}$ denoting the restriction of a map $A$ from two spaces $V$ and $W$ to subspaces $V_0\subset V, W_0\subset W$, with inverse (if exists) denoted by $A_{V_0\downarrow W_0}^{-1}$, shorthand for $(A_{V_0\downarrow W_0})^{-1}$. 
\begin{proposition}\label[proposition]{prop:prelimOverdeter}Assume the map $\cF$ from $\Omega\subset T_0$ to $\EE$ and the projection function $\Phi:\Omega\mapsto\Lin(\EE, \EE)$ are smooth. Assume for $\xi\in \Omega$
      \begin{equation}\label{eq:cFsection}
    \Phi(\xi)\cF(\xi) = \cF(\xi).
\end{equation}
      Set $\EE_{\xi}:=\Imag(\Phi(\xi))\subset \EE$. If $\Phi(\xi)\cF'(\xi)$ is a bijection from $T_0$ to $\EE_{\xi}$, for $\omega\in \EE$, let $\cF^{\natural}(\xi)\omega\in T_0$ be the unique solution $\eta$ of $\Phi(\xi)\cF'(\xi)\eta = \Phi(\xi)\omega$, thus 
$\cF^{\natural}(\xi)_{\EE_{\xi}\downarrow T_0} = (\Phi(\xi)\cF'(\xi))_{T_0\downarrow \EE_{\xi}}^{-1}$ and     
  \begin{equation}\label{eq:Fnatural}
    \cF^{\natural}(\xi) := (\Phi(\xi)\cF'(\xi))_{T_0\downarrow \EE_{\xi}}^{-1}\Phi(\xi)\in \Lin(\EE, T_0).
  \end{equation}
  1) For each $\xi$, the linear map $\cF^{\natural}(\xi)\in \Lin(\EE, T_0)$, if exists, is a left inverse of $\cF'(\xi)$ and
  \begin{gather}\cF^{\natural}(\xi)\cF'(\xi) = I_{T_0}\\
    \Phi(\xi)\cF'(\xi)\cF^{\natural}(\xi) = \Phi(\xi).
  \end{gather}
Assume $\cF^{\natural}(\xi)$ is defined in a domain $D_1\subset \Omega$, then it is differentiable and for $\xi, \eta\in T_0$ and $\omega \in \EE$
\begin{equation}  
\label{eq:JFnat}
    (\cF^{\natural})'(\xi; \eta)\omega =
 \cF^{\natural}(\xi)\{-\cF^{(2)}(\xi;\eta, \cF^{\natural}(\xi)\omega) +
 \Phi'(\xi;\eta)[\omega -  \cF'(\xi) \cF^{\natural}(\xi)\omega]\}.
\end{equation} 
  2) Assume $\cF^{\natural}(\xi)$ is defined in a domain $D_1\subset \Omega$. Then for $(\xi, Y)\in D_1\times \EE$
\begin{equation}\label{eq:Fdiamond}
  \Phi(\xi)Y - \cF(\xi)=0 \Leftrightarrow  \cF^{\natural}(\xi)(Y - \cF(\xi)) = 0.
  \end{equation}
For $\xi_f\in \EE$, there is an implicit-function solution $\xi=\cF_{\diamond}(Y)$ for the second equation near $(\xi_f, \cF(\xi_f))\in D_1\times \EE$, with $Y\in\EE$, $\xi\in T_0$, where $\cF_{\diamond}$ is of class $C^k$ if $\cF$ and $\Phi$ are. If $(\xi, Y)$ satisfies \cref{eq:Fdiamond} then $(\xi, Z)$ also satisfies \cref{eq:Fdiamond} for $Z = \Phi(\xi)Y$. For $\omega\in \EE$, we have
\begin{gather}
  (I_{T_0} - \cF^{\natural}(\xi)\Phi'(\xi; .; Y))\cF_{\diamond}'(Y)\omega :=
\cF_{\diamond}'(Y)\omega - \cF^{\natural}(\xi)\Phi'(\xi; \cF_{\diamond}'(Y)\omega) Y   = \cF^{\natural}(\xi)\omega\label{eq:bJFdiam}\\
  \cF_{\diamond}'(Z)\omega = \cF^{\natural}(\xi)\omega.\label{eq:bJFdiamZ}
 \end{gather}
Assume $\cF$ and $\Phi$ are smooth, then for $\omega\in\EE$
\begin{equation}
  \cF_{\diamond}^{(2)}(Z; \omega^{[2]}) = \cF^{\natural}(\xi)\{-\cF^{(2)}(\xi; (\cF^{\natural}(\xi)\omega)^{[2]}) +2\Phi'(\xi, \cF^{\natural}(\xi)\omega)\omega+ \Phi^{(2)}(\xi; (\cF^{\natural}(\xi)\omega)^{[2]})Z\}\label{eq:jFDia2}
\end{equation}  
3) If $(I_{T_0} - \cF^{\natural}(\xi)\Phi'(\xi; .; Y))^{-1}$ exists and is bounded in a region $D_1$, and $\cF$ and $\Phi$ are of class $C^k$, $k\geq 1$ with bounded derivatives, then $\cF_{\diamond}$ is also of class $C^k$ with bounded derivatives.
\end{proposition}  

We now formulate the convergence theorem with a nontrivial retraction. We will see the retraction in \cref{prop:fSretract} will appear in RQIs. For the case of the sphere \cref{eq:fSsphere}, it has the form $\fS_{s}:(\xi, \delta) \mapsto \xi + \frac{1+|\xi|^2}{1 + \xi^{\ft}(\xi-\delta)}\delta$, for $\xi, \delta\in T_0$, and is only defined for $\delta$ satisfying \cref{eq:condrho12}. Thus, we need to address retractions of this type.
\begin{theorem}\label{thm:mainOver}Assume the function $\cF$, the projection $\Phi$ are smooth satisfying $\Phi(\xi)\cF(\xi) =\cF(\xi)$ for $\xi\in\Omega\subset T_0$ as in \cref{eq:cFsection}.\hfill\break
  1) Assume $(\Phi(\xi)\cF'(\xi))_{T_0\downarrow \EE_{\xi}}^{-1}$ is defined in a ball $B(\xi_*, \rho)$ in $T_0$, with $\cF(\xi_*) = 0$. For $\rho, \rho_2 >0$, assume the map $\fS$ from $D_1 = B(\xi_*, \rho)\times B(0, \rho_2)\subset T_0\times T_0$ to $T_0$ of class $C^2$ has the Taylor expansion in $\delta$
  \begin{equation}\fS(\xi, \delta) = \xi + \delta + \fS_2(\xi, \delta)
  \end{equation}
  such that there is a constant $c_{\fS, 2}$ with $|\fS_2(\xi, \delta)|\leq c_2|\delta|^2$ in $D_1$. Then there exists a radius $\rho_3$ such that if the starting point $\xi_0$ satisfies $|\xi_0-\xi_*| < \rho_3$ then the iteration defined by
  \begin{gather}\xi_{i+1} =\fS(\xi_i, \delta_{2, i}) \text{ for }
    \delta_{2, i} := -\cF^{\natural}(\xi_i)\cF(\xi_i)
  \end{gather}
is well-defined and converges quadratically to $\xi_*$.\hfill\break
  2) If $\fS$ is of class $C^3$ and the 3-term Taylor series of $\fS$ in $\delta$ is given by \begin{equation}\fS(\xi, \delta) = \xi + \delta + \frac{1}{2}\fS(\xi, .)^{(2)}(0; \delta^{[2]}) +  \fS_3(\xi, \delta)
\end{equation}
with $|\fS_3(\xi, \delta)| < c_{\fS,3}|\delta|^3$ and $\cF$ is of class $C^3$. If $\Phi(\xi_*)\cG(\xi_*, \delta^{[2]})=0$ for all $\delta\in T_0$ where
\begin{equation}\label{eq:bGstar}\cG(\xi, \delta^{[2]}) := \cF^{(2)}(\xi; \delta^{[2]}) + \cF'(\xi; \fS(\xi, .)^{(2)}(0; \delta^{[2]}))\end{equation}
then the Newton iteration in 1) converges cubically to $\xi_*$.  \hfill\break
3) With the smoothness assumption in 2), but we do not assume $\Phi(\xi_*)\cG(\xi_*, \delta^{[2]})=0$, then there exists a radius $\rho_4$ such that for $|\xi_0-\xi_*| < \rho_4$, the following iteration converges cubically to $\xi_*$
\begin{gather}\xi_{i+1} = \fS(\xi_i, \delta_{3,i})\\
  \delta_{3,i} := \delta_{2, i} - \frac{1}{2}\cF^{\natural}(\xi)\cG(\xi_i, \delta_{2, i}^{[2]})
\end{gather}
\end{theorem}
The proof is in \cref{appx:prf:mainOver}.

\section{Constrained equations and Rayleigh\slash Rayleigh-Chebyshev iterations}\label{sec:RQIRC}
We now consider a feasible set $\cM$ defined by a full range constraint $\bC(\bx) = 0$ with a retraction $\bR$. Let $\Pi$ be an affine projection function from $\cM$ to $\Lin(\EE, \EE)$, $\Pi(\bx)^2=\Pi(\bx)$ for $\bx\in\cM$. Denote $\EE_{\bx}=\Imag(\Pi(\bx)), T_{\bx}=\Null(\bC'(\bx))$. In the subset of $\cM$ considered, we assume
\begin{equation} \label{eq:rankCons}\dim\EE_{\bx} = \dim(\Imag(\Pi(\bx))) = \dim T_{\bx} =\dim(\EE) - \dim(\EL).\end{equation}
We try to solve $\bF(\bx) = 0$ for a function $\bF$ from $\cM$ to $\EE$ satisfying the condition
\begin{equation}\label{eq:secCons}\Pi(\bx)\bF(\bx) = \bF(\bx). \end{equation}
If $(\Pi(\bx)\bF'(\bx))_{T_{\bx}\downarrow E_{\bx}}$ is invertible, (i.e. the system in \cref{claim:cl1} is nondegenerate), define
\begin{equation}\bF^{\natural}(\bx) := (\Pi(\bx)\bF'(\bx))_{T_{\bx}\downarrow E_{\bx}}^{-1}\Pi(\bx).\end{equation}
\begin{proposition}\label[proposition]{prop:Fdef0} Fix a feasible point $\bx_f\in \cM$, let $\bR$ be a retraction around $\bx_f$ with a retraction cap $\bC_{\bx_f}(\rho)\subset \cM$, parametrized by a retraction ball $B_{\bx_f}(\rho)\subset T_{\bx_f}$ via $\tK_f :=\bR(\bx_f, .)$. Let $\bF$ be a map from $\cM$ to $\EE$, $\Pi$ be an affine projection. For $\bx=\bR(\bx_f, \xi)\in C(\bx_f, \rho)$, consider the Newton increment and step
    \begin{gather}
    \bx^{\cN_T} = -\bF^{\natural}(\bx)\bF(\bx)= -(\Pi(\bx)\bF'(\bx))_{T_{\bx}\downarrow E_{\bx}}^{-1}\Pi(\bx)\bF(\bx)\in T_{\bx}\label{eq:XNewton01}\\
    \bx^{\cN} = \bR(\bx, \bx^{\cN_T}) \in \cM\label{eq:XNewton02}
    \end{gather}
assume they are well-defined and $\bx^{\cN}\in C_{\bx_f}(\rho)$. For $\xi\in B_{\bx_f}(\rho)$, set $\Phi(\xi) := \Pi(\bR(\bx_f, \xi))$, $\cF(\xi) := \bF(\bR(\bx_f, \xi))$. Then the Newton-type iteration $\xi^{\cN_f}$ of the problem $\cF(\xi) = 0$
    \begin{gather}
      \xi^{\cN_f} := \fS(\xi, \xi^{\Delta})\in T_{\bx_f}\label{eq:XNewton04}\\      
      \xi^{\Delta} := - \cF^{\natural}(\xi)\cF(\xi)\in T_{\bx_f}\label{eq:XNewton05}
    \end{gather}
with $\fS$ defined in \cref{eq:FS} corresponds to the Newton iteration $\bx^{\cN}$. That means
\begin{equation}    
      \bx^{\cN} =\bR(\bx_f, \xi^{\cN_f})\label{eq:XNewton03}.\\
\end{equation}    
\end{proposition}
\begin{proof}From the chain rule, $\cF'(\xi)$ and its left-inverse $\cF^{\natural}(\xi)$ in \cref{eq:Fnatural} are
  \begin{gather}
    \cF'(\xi) = \bF'(\bR(\bx_f, \xi))\bR(\bx_f, .)'(\xi)\in\Lin(T_{\bx_f}, \EE)\\
    \cF^{\natural}(\xi) =
    (\bR(\bx_f, .)'(\xi))^{-1}(\Phi(\xi)\bF'(\bR(\bx_f, \xi)))^{-1}_{T_{\bx_f}\downarrow E_{\bx}}\Phi(\xi)
= (\tK_f'(\xi))^{-1}\bF^{\natural}(\bx)
    \in \Lin(\EE, T_{\bx_f})\label{eq:bFpNat}
  \end{gather}
where $(\bR(\bx_f, .)'(\xi))$ moves to the front in the inverse. For \cref{eq:XNewton03,eq:XNewton04,eq:XNewton05}, with $X = \bR(\bx_f, \xi)$
  $$\begin{gathered}
\tK_f(\xi^{\cN_f}) =  \tK_f(\fS(\xi, \xi^{\Delta})) = \tK_f\circ \tK_f^{-1}(\bR(\tK_f, \xi), \tK_f'(\xi;\xi^{\Delta}))\\
    = \bR(\tK_f, \xi), -\tK_f'(\xi)[(\Pi(\tK_f(\xi))\bF'(\tK_f(\xi))\tK_f'(\xi))^{-1}_{T_{\bx_f}\downarrow E_{\bx}}\Pi(\tK_f(\xi))\bF(\tK_f(\xi))]
    =\bR(\tK_f(\xi), \bx^{\cN_T} ) =\bx^{\cN}.\qedhere
  \end{gathered}$$
\end{proof}

We now translate the theorems in \cref{sec:HOGaussNewton} to results on constrained iterations.
\begin{theorem}[Local convergence]\label{thm:LocalConstrained}Assume $\bF$ and $\Pi$ satisfy \cref{eq:rankCons,eq:secCons}. Assume $\bF(\bx_*)=0$ and $(\Pi(\bx_*)\bF'(\bx_*))_{T_{\bx_*}\downarrow E_{\bx_*}}$ is invertible. Let $\bR$ be a retraction around $\bx_*$. Then there is a radius $\rho >0$ defining a retraction ball $B_{\bx_*}(\rho)\subset T_{\bx_*}$ and cap $C_{\bx_*}(\rho)\subset \cM$ such that for $\bx_0\in C_{\bx_*}(\rho)$, the iteration $\bx_{i+1} = \bx_{i}^{\cN}$ defined in \cref{prop:Fdef0} is well-defined and converges to $\bx_*$.

If $\Pi(\bx_*)\bG(\bx_*; \eta^{[2]}) = 0$ for $\eta\in T_{\bx_*}$, where $\bG$ is a tensor-valued function on $\cM$ defined by
  \begin{equation}\label{eq:GRQI}
    \bG(\bx; \eta^{[2]}) = \bF^{(2)}(\bx, \eta^{[2]}) + \bF'(\bx; \bR(\bx, .)^{(2)}(0;\eta^{[2]}))\quad \text{ for }\eta \in T_X
  \end{equation}
then the Newton iteration $\bx_{i+1} = \bx_{i}^{\cN}$ converges cubically. Define the Chebyshev increment
  \begin{equation}
    \bx^{\cN_C} = \bx^{\cN_T} -\frac{1}{2}(\Pi(\bx)\bF'(\bx))_{T_{\bx}\downarrow \EE_{\bx}}^{-1}\Pi(\bx)\bG(\bx; \bx^{\cN_T})\in T_{\bx}
  \end{equation}
Then there is $\rho_C > 0$ such that for $\bx_0\in C_{\bx_*}(\rho_C)$, the Rayleigh-Chebyshev iteration $\bx_{i+1} = \bR(\bx_i, \bx_i^{\cN_C})$ converge cubically to $\bx_*$.
\end{theorem}
\begin{proof}Using \cref{prop:Fdef0}, the Newton case follows from \cref{thm:mainOver}. Set $\cF(\xi) = \bF(\bR(\bx_f, \xi))$ for a feasible point $\bx_f$, then differentiate $\cF'(\xi)\delta = \bF'(\bR(\bx_f, \xi))\bR(\bx_f, .)'(\xi; \delta)$ in direction $\delta$
  $$\begin{gathered}\cG(\xi, \delta^{[2]}) = \cF^{(2)}(\xi; \delta^{[2]}) + \cF'(\xi; \fS(\xi, .)^{(2)}(0; \delta^{[2]})) = \bF^{(2)}(\bR(\bx_f, \xi); [\bR(\bx_f, .)'(\xi; \delta)]^{(2)}) + \\
\bF'(\bR(\bx_f, \xi))\bR(\bx_f, .)^{(2)}(\xi; \delta^{[2]}) 
+ \bF'(\bR(\bx_f, \xi))[\bR(\bR(\bx_f, \xi), . )^{(2)}(0; [\bR(\bx_f, .)'(\xi; \delta)]^{[2]}) -    \bR(\bx_f, .)^{(2)}(\xi; [\delta]^{[2]})]\\
=\bF^{(2)}(\bx; \eta^{(2)}) + \bF'(\bx)\bR(\bx, .)^{(2)}(0; \eta^{[2]}) =\bG(\bx, \eta^{[2]})
  \end{gathered}$$
  if $\bx = \bR(\bx_f, \xi), \eta  = \bR(\bx_f, .)'(\xi; \delta)$ and $\cG$ defined in \cref{eq:bGstar}. Cubic convergence of the Newton algorithm if $\bG$ vanishes follows. The statement for Rayleigh-Chebyshev also follows.
\end{proof}
\subsection{Rayleigh quotient iteration for constrained equations}\label{subsec:RQItheo}
As mentioned, for the equation \cref{eq:LMR0}, if $\bLlbd(\bx, \blbd)$ is onto, as a linear map from $\EL$ to $\EE$, then there is a left-inverse $\bL^-_{\lambda}(\bx, \blbd)$ of $\bLlbd(\bx, \blbd)$. If it is constructed smoothly, then with the Rayleigh quotient $\lambda = \cR(\bx)$ solving $\bL^-_{\lambda}(\bx, \blbd)\bL(\bx, \blbd)=0$, the projection $\Pi(\bx):=I_{\EE} -
\bL_{\lambda}(\bx, \cR(\bx))\bL_{\lambda}^-(\bx, \cR(\bx))$ satisfies
\begin{equation}\label{eq:Rchoice}
  \begin{gathered}
    \Pi(\bx)\bL(\bx, \blbd) = \bL(\bx, \blbd),\\
    \Pi(\bx)\bLlbd(\bx, \blbd) = 0.
\end{gathered}\end{equation}
Let $\bF(\bx) = \bL(\bx, \cR(\bx))$, then $\Pi(\bx)\bF(\bx) = \bF(\bx)$, thus, the results of the previous section apply. Using $\Pi(\bx)\bLlbd(\bx, \cR(\bx)) = 0$ to simplify $\Pi(\bx)\bF'(\bx)$ and $\Pi(\bx)\bG(\bx)$, the following theorem follows from the chain rule. The second derivatives of $\bL$ are bilinear maps, $\bL_{\bx\bx}$ operates on two copies of $\EE$ via the variable $\eta$ below, $\bL_{\bx\blbd}$ operates on $\EE\times\EL$ via $\eta\in\EE$ and $\cR'(\bx; \eta)\in\EL$, and $\bL_{\blbd\blbd}$ on $\cR'(\bx; \eta)^{[2]}\in\EL^2$.
\begin{theorem}\label{thm:RQIGEN} If the Rayleigh quotient $\lambda =\cR(\bx)$ and the projection $\Pi$ satisfies \cref{eq:Rchoice}, with $\bL$ and $\cR$, $\bC$, $\bR$, $\Pi$ are smooth, then with $\lambda =\cR(\bx)$, the Newton increment $\bx^{\cN_T}\in T_{\bx}$ of the RQI iteration $\bx_{i+1} = \bR(\bx_i, \bx^{\cN_T})$ is the solution to the equation
  \begin{equation}\label{eq:RQILambda}
   \Pi(\bx)\bLx(\bx, \lambda)\bx^{\cN_T}  = -\bL(\bx, \lambda)
  \end{equation}
  and the convergence criteria for $\bF = \bL(\bx, \cR(\bx))=0$ in \cref{thm:LocalConstrained} apply. Let
  \begin{equation}\bF^{\natural}(\bx)\omega = (\Pi(\bx)\bLx(\bx, \lambda))_{T_{\bx}\downarrow E_{\bx}}^{-1}\Pi(\bx)\omega\quad \text{ for }\omega\in \EE.
    \end{equation}
  Assuming $\bC$, $\bR$, $\Pi$, $\bL, \cR$ are smooth. For $\eta\in T_{\bx}$, set $\nu = \bR(\bx, .)^{(2)}(0;\eta^{[2]})$ and
\begin{equation}\label{eq:BGRQI}\begin{gathered} \bG_{\bL}(\bx; \eta^{[2]}) = \bL_{\bx\bx}(\bx, \lambda)\eta^{[2]} + 2\bL_{\bx\blbd}(\bx, \lambda)[\eta, \cR'(\bx; \eta)]+
    \bL_{\blbd\blbd}(\bx, \blbd)\cR'(\bx; \eta)^{[2]}  + \bLx(\bx, \lambda)\nu
  \end{gathered}
\end{equation}
then $\Pi(\bx)\bG_{\bL}(\bx; \eta^{[2]})=0$ if and only if $\Pi(\bx)\bG(\bx; \eta^{[2]})=0$ in \cref{eq:GRQI}. Thus, if $\Pi(\bx_*)\bG_{\bL}(\bx_*; \eta^{[2]})=0$ for all $\eta\in T_{\bx_*}$ and $\bF^{\natural}(\bx_*)$ exists for a solution $\bx_*$ of $\bF(\bx) = 0$, the RQI above converges cubically. If $\bF^{\natural}(\bx_*)$ exists then the Rayleigh-Chebyshev iteration $\bx_{i+1} = \bR(\bx_i, \bx^{\cN_C})$ with
  \begin{equation}
    \bx^{\cN_C} = \bx^{\cN_T} -\frac{1}{2}\bx^{\tau} = \bx^{\cN_T} -\frac{1}{2}(\Pi(\bx)\bLx(\bx))_{T_{\bx}\downarrow \EE_{\bx}}^{-1}\Pi(\bx)\bG_{\bL}(\bx; \bx^{\cN_T})\in T_{\bx}
  \end{equation}
 converges cubically to $\bx_*$ if $\bx_0\in C_{\bx_*}(\rho)$ for some $\rho> 0$. If $\bLx(\bx, \lambda)\in \Lin(\EE, \EE)$ is invertible and $\bC'(\bx)\bLx(\bx, \lambda)^{-1}\bLlbd(\bx, \lambda)$ is invertible, $\bx^{\cN_T}$ and $\bx^{\tau}$ could be computed in Schur form
  \begin{gather}
    \bx^{\cN_T} =-\cS(\bLx(\bx, \lambda)^{-1}\bL(\bx,\lambda)) \\
    \bx^{\tau} = \cS(\bLx(\bx, \lambda)^{-1}\bG_{\bL}(\bx))\\
    \cS\omega := \omega - \bLx(\bx, \lambda)^{-1}\bLlbd(\bx, \lambda)
(\bC'(\bx)\bLx(\bx, \lambda)^{-1}\bLlbd(\bx, \lambda))^{-1}\bC'(\bx)\omega\quad\quad\text{ for }\omega\in \EE.
  \end{gather}
\end{theorem}
If $\bL$ is affine in $\blbd$, $\bL(\bx, \blbd) = A(\bx) - H(\bx)\lambda$ for two functions $A(\bx)$ and $H(\bx)$, then $\bLx(\bx, \blbd) = A'(\bx) - H'(\bx)\lambda$. Dropping the variable name $\bx$ for brevity, if $H^-$ is a left inverse of $H$, we can take $\bL^-_{\blbd}$ to be $-H^-$, giving us a Rayleigh quotient $\lambda=H^-A$ and $\bF = A - HH^-A = \Pi A$ with $\Pi = I_E - HH^-$. Thus, in \cref{lem:cl2}, we can use $\bA$ instead of $\bF$ in the Schur form solution
\begin{equation}\label{eq:schurAH}\bx^{\cN_T} = -\bLx^{-1}A +\bLx^{-1}H(\bC'\bLx^{-1} H)^{-1}\bC'\bLx^{-1}A.
\end{equation}
\section{Cubic convergence and second covariant derivative}\label{sec:geomcubic}
In \cref{sec:geometry}, with the projection function $\Pi$, we defined the vector bundle $E_{\Pi}=\cup_{\bx\in\cM}\Imag(\Pi(\bx))$ with a {\it connection} $(\nabla_{\tc}\bF)(\bx) = \Pi(\bx)\bF'(\bx, \tc(\bx))\in T_X$ for a {\it section} $\bF$, $\Pi(\bx)\bF(\bx) = \bF(\bx)$ and a vector field $\tc$, $\tc(\bx)\in T_{\bx}$. In the unconstrained case, the Chebyshev iteration
$$\bx_{i+1} = \bx_i - \bF'(\bx)^{-1}\bF(\bx) -\frac{1}{2}\bF^{(2)}(\bx; [\bF'(\bx)^{-1}\bF(\bx)]^{[2]})$$
converge cubically, as seen in \cref{sec:HOGaussNewton}. It is natural to ask if there is a relationship between the Rayleigh-Chebyshev term and $\nabla$ in the constrained case. The answer is yes, with the concept of a {\it second covariant derivative}.

If $\bC$ is the constrained function for $\cM$, for two vector fields $\tc_1, \tc_2$, functions from $\cM$ to $\EE$ satisfying $\bC'(\bx)\tc_i(\bx) = 0$, $i=1,2,\bx\in \cM$, then $\tc_2'(\bx; \tc_1(\bx))$ is not a vector field, as taking derivative of $\bC'(\bx)\tc_2(\bx)=0$ in direction $\tc_1(\bx)$ does not give $\bC'(\bx)\tc_2'(\bx; \tc_1(\bx))=0$ but
\begin{equation}\bC'(\bx)\tc_2'(\bx; \tc_1(\bx)) + \bC^{(2)}(\bx; \tc_1(\bx), \tc_2(\bx))=0.\label{eq:Cc1c2}\end{equation}
An adjustment of the form $\tc_2'(\bx, \tc_1(\bx)) + \Gamma(\bx, \tc_1(\bx), \tc_2(\bx))$ is needed to get a vector field, called a {\it covariant derivative} of $\tc_2$ at $\bx$ along $\tc_1$. As seen in \cref{sec:geometry}, a projection function to $T\cM$ defines a covariant derivative (connection). A retraction $\bR$ also defines a covariant derivative $\nabla^{\bR}$:

\begin{proposition}\label[proposition]{prop:retractConnect}For a retraction $\bR$ and two vector fields $\tc_1, \tc_2$, $\nabla_{\tc_1}^{\bR}\tc_2$ defined by
  \begin{equation}  \nabla_{\tc_1}^{\bR}\tc_2:\bx \mapsto\tc_2'(\bx; \tc_1(\bx)) - \bR(\bx, .)^{(2)}(0; \tc_1(\bx), \tc_2(\bx))
  \end{equation}
is a vector field, thus $\nabla^{\bR}$ is a covariant derivative.
\end{proposition}  
\begin{proof}
Differentiate $\bC'(\bR(\bx, t\tc_1(\bx)))\bR(\bx, .)'(t\tc_1(\bx); \tc_2(\bx)) = 0$ (from \cref{lem:vectorTrans}) in $t$ then set $t=0$ we get
  $$\begin{gathered}
    \bC^{(2)}(\bx; \tc_1(\bx), \tc_2(\bx)) + \bC'(\bx)\bR(\bx, .)^{[2]}(0, \tc_1(\bx), \tc_2(\bx)) =0.
  \end{gathered}$$
  Subtracting from \cref{eq:Cc1c2}, we get
  $$\bC'(\bx)\{\tc_2'(\bx; \tc_1(\bx)) - \bR(\bx, .)^{(2)}(0; \tc_1(\bx), \tc_2(\bx))\} =0.$$
This shows $\nabla_{\tc_1}^{\bR}\tc_2(\bx)$ is tangent to $\cM$ at $\bx$.
\end{proof}

The association of a connection with a retraction is discussed in \cite[section 8.1.1]{AMS_book}. Here, we give an explicit formula. A naive attempt for a geometric expression of the Chebyshev term is to consider an adjustment like $\frac{1}{2}\nabla_{\eta}(\nabla_{\eta}\bF)$ , with $\nabla$ is the connection from $\Pi$ in \cref{sec:geometry} and $\eta$ is the Newton increment. To evaluate the outer
$\nabla_{\eta}$, we need the inner to be a section, for this we need to replace $\eta$ with a vector field. For two vector fields $\tc_1, \tc_2$
$$\begin{gathered}(\nabla_{\tc_1}(\nabla_{\tc_2}\bF))(\bx) = \Pi(\bx)\{\bx\mapsto \Pi(\bx)\bF'(\bx; \tc_2(\bx))\}'(\bx; \tc_1(\bx))\\
=\Pi(\bx)\{\Pi'(\bx; \tc_1(\bx))\bF'(\bx)\tc_2(\bx) + \Pi(\bx)\bF^{(2)}(\bx; \tc_1(\bx), \tc_2(\bx)) + \Pi(\bx)\bF'(\bx; \tc_2'(\bx; \tc_1(\bx)))\}.
\end{gathered}$$
The last term involves $\tc_2'(\bx;\tc_1(\bx))$, hence dependent on the extensions of $\eta$ to vector fields. But $\nabla_{\nabla_{\tc_1}^{\bR}\tc_2}\bF(\bx)=\Pi(\bx)\bF'(\bx; \nabla_{\tc_1}^{\bR}(\tc_2\bx))$ is well-defined, hence $\nabla^2_{\tc_1, \tc_2}\bF(\bx)$ given below is not dependent on derivatives of $\tc_1$ or $\tc_2$, it is the {\it second covariant derivative} using $\nabla$ and $\nabla^{\bR}$
\begin{equation}\begin{gathered}\nabla^2_{\tc_1, \tc_2}\bF(\bx) :=
(\nabla_{\tc_1}(\nabla_{\tc_2}\bF))(\bx) - \nabla_{\nabla_{\tc_1}^{\bR}\tc_2}\bF(\bx)\\
=\Pi(\bx)\{\Pi'(\bx, \tc_1(\bx))\bF'(\bx)\tc_2(\bx) + \bF^{(2)}(\bx; \tc_1(\bx), \tc_2(\bx)) + \bF'(\bx)\bR(0, .)^{[2]}(0; \tc_1(\bx), \tc_2(\bx))\}.
\end{gathered}\end{equation}
Therefore, we can define $\frac{1}{2}\nabla^2_{\eta, \eta}\bF(\bx)$ as above, for vector fields satisfying $\tc_1(\bx) = \tc_2(\bx) = \eta$. \cite{castro} uses this expression as the third-order adjustment for the geodesic retraction. Compare with \cref{eq:GRQI}, we have an extra term $\Pi(\bx)\Pi'(\bx, \eta)\bF'(\bx)\eta$. Differentiate $\Pi(Y)\Pi'(Y; \tc(Y))\bF(Y)=0$ (\cref{eq:secondFForm}) in direction $\tc(Y)$, $Y\in\cM$
$$\begin{gathered}\Pi'(Y; \tc(Y))\Pi'(Y; \tc(Y))\bF(Y) +\Pi(Y)\Pi^{(2)}(Y; \tc(Y)^{[2]})\bF(Y)
+ \\ \Pi(Y)\Pi'(Y; \tc'(Y)\tc(Y))\bF(Y) 
 + \Pi(Y)\Pi'(Y; \tc(Y))\bF'(Y)\tc(Y) = 0.
\end{gathered}$$
At a solution $\bx_*$, $\bF(\bx_*) = 0$, the first three terms are zero, hence $\Pi(\bx_*)\Pi'(\bx_*, \eta)\bF'(\bx_*)\eta=0$ for any tangent vector $\eta$ at $\bx_*$ if $\tc(\bx_*)=\eta$. Thus, this term does not affect the order of convergence.
\section{Applications}\label{sec:examples}
The codes for these sections are found in folder \href{https://github.com/dnguyend/rayleigh_newton/tree/master/colab}{colab} in \cite{rayleigh_newton_code}, implemented in Python and Julia (for two examples showing cubic convergence of Rayleigh-Chebyshev iterations). 
\subsection{Complex tensor eigenpairs}\label{subsec:tensor}
For our purpose, a tensor is a vector-valued homogeneous function. An example of the real tensor eigenpair problem is finding critical points of a scalar homogeneous function $\hcT$ under the constraint $\hcB=1$. Here, $\hcB$ is a scalar homogeneous polynomial, which leads to the equation $\grad\hcT = \lambda \grad\hcB$. We will focus on equations of the form $\cT = \lambda \cB$ explained below, the constraints imposed may be unrelated to $\cB$.

We consider a complex-valued homogeneous vector function $\cT$, which may not be a gradient of a scalar function. Each entry of $\cT$ is a scalar homogeneous polynomial of order $m-1$. The coefficients $(t_{i_1\cdots i_m})_{1\leq i_j\leq n, 1\leq j\leq m}$ of entry $\cT_{i_1}, i_1=1\cdots n$ can be arranged to what is called a multidimensional array (also called a tensor), we call $\cT$ a tensor of dimension $n$ order $m$. For $k$ vectors $X_{n-k+1}\cdots X_n$ of size $n$ in the base field, contraction means taking sum over the last $k$ indices of $(t_{i_1\cdots i_m})$ multiply by entries of $X_j$'s, $j=n-k+1\cdots n$, the uncontracted entries are written $I$. Write $\bx^{[k]}$ for $\bx$ repeated $k$ times. Abuse of notation, $\cT(\bx^{[m-1]})=\cT(I, \bx^{[m-1]})$ is the vector homogeneous function evaluated at $\bx$, its Jacobian is $(m-1)\cT(I, I, \bx^{[m-2]})$, while $\hcT(\bx^{[m]})$ is the scalar function in the symmetric case, where $\grad\hcT = m \cT$ (We mostly follows \cite{JaffeWeissNadler}).

Let $\cB$ be a homogeneous polynomial vector function of order $d-1$, the tensor eigenpair problem in \cref{eq:tensordef} solves for pairs $(\bx, \lambda)$ satisfying $\cT(Z^{[m-1]})=\blbd \cB(Z^{[d-1]})$.

We will focus on the case $d=2$, much of what is discussed here holds for $d\geq 2$. The number of complex eigenpairs is given as by $\sum_{i=0}^{n-1}(m-1)^i(d-1)^{n-i-1}$ (\cite{Homotopy,CARTWRIGHT2013942}), with certain counting convention described there. An important case is the Z-pairs \cite{LiqunQi}, \cite{Lim}, where $\cB(Z)=Z$
  \begin{equation}\bL(Z, \blbd) = \cT(Z^{[m-1]})-\blbd Z = 0\label{eq:tensoreq}
  \end{equation}
  with the normalization $Z^{\ft}Z=1$ in the real case \cite{LiqunQi}. For the complex case, we normalize by $Z^*Z=1$ ($*$ is the Hermitian transpose), and call the pairs unitary Z-pairs (UZ-pairs).
  
  There are a few methods to find all real eigenpairs \cite{Homotopy,CuiDaiNie}, with the first citation using the homotopy method is state-of-the-art, which also finds complex pairs. In \cite{JaffeWeissNadler}, the authors proposed a RQI-type algorithm, called $\ONCM$, which empirically finds all real Z-pairs in less time than the homotopy method, however, there is no certification that all real pairs are found, as there is no formula to count real pairs. Using the RQI method described in this paper, we extend the $\ONCM$ method to the complex case, using the count $\sum(m-1)^i$ for complex (Z-pairs) to decide if we have found all eigenpairs, if all pairs are distinct, up to scaling.

  If $m >2$, we can assume $\blbd$ is real, or $\EL=\R$, by scaling $\lambda$ by $\exp(-i\theta)$ and $Z$ by $\exp(-i\theta/(m-2))$, with $\theta$ is the angle of $\lambda$ in polar coordinates. We treat $Z$ as a real vector of dimension $2n$.
From $\bC(Z) = Z^*Z-1$, we have $\bC'(Z)\eta= 2\real(Z^*\eta)$, thus a tangent vector satisfies $\real(Z^*\eta) = 0$. We have $\bLlbd(Z, \lambda)(\delta) = -Z\delta$, $\delta\in\R$. A left inverse is $\bLlbd^-:\omega\mapsto -\real(Z^*\omega)\in \EL, \omega\in\EE$, thus
$$\Pi^-(Z)\omega = \bLlbd(\bx, \blbd)\bLlbd^-(\bx, \blbd)\omega= Z\real(Z^* \omega)\quad\quad\text{ for } \omega\in \EE.$$
To use \cref{thm:RQIGEN}, we solve for the Rayleigh quotient $\lambda=\cR(Z)$ from $\bLlbd^-(\bx, \lambda)\bL(Z, \lambda)=0$, or $\real(Z^*\cT(Z^{[m-1]}) - \lambda \real Z^*Z) = 0$
\begin{gather}\cR(Z) = \real(Z^*\cT(Z^{[m-1]}))
\end{gather}
The projection is $\Pi(Z)\omega = \omega - Z\real(Z^*\omega)$ for $\omega\in \EE$. The Newton increment $Z^{\cN_T}$ satisfies
\begin{equation}
\Pi(Z)\bL_{Z}(Z,\blbd)Z^{\cN_T}=  \Pi(Z)((m-1)\cT(\bI, \bI, Z^{[m-2]}) - \blbd\bI)Z^{\cN_T} = - \Pi(Z)\cT(Z^{[m-1]})
\end{equation}
The Newton form \cref{eq:NewtonForm} is the complex version of the algorithm $\ONCM$ op. cit. The tangent space is the same as the range $E_{Z}$ of $\Pi(Z)$, both are nullspace of $\omega\mapsto\real Z^*\omega$. In \cref{eq:schurAH}, with $A=\cT(Z^{[m-1]})$, $H=Z$, set $\zeta = \bLx^{-1}H, \nu = \bLx^{-1}A$, the Schur form solution is
\begin{equation}\begin{gathered}  Z^{\cN_T} =  -\nu + \zeta\real(Z^*\zeta))^{-1} \real(Z^*\nu)  \\
    \text{ where }    \begin{bmatrix}\zeta & \nu\end{bmatrix} =
 ((m-1)\cT(\bI, \bI, Z^{[m-2]}) - \blbd\bI)^{-1}\begin{bmatrix}Z & \cT(Z^{[m-1]})\end{bmatrix}\in \C^{n\times 2}.
\end{gathered}      
\end{equation}
\begin{algorithm}
  \begin{algorithmic}[1]
    \State{Input: $Z_0$ with $|Z_0| = 1$}
\For{$i=0,1,\cdots$}    
\State Compute $\blbd_{i} = \cR(Z_i) = \real(Z_i^*\cT(Z_i^{[m-1]}))$\;
\State Solve $((m-1)\cT(\bI, \bI, Z_i^{[m-2]}) - \blbd_i\bI) \begin{bmatrix}\zeta & \nu\end{bmatrix} = \begin{bmatrix}Z_i & \cT (Z_i^{[m-1]})\end{bmatrix}$\;
\State Compute $\eta \gets -\nu + \zeta\real(Z_i^*\zeta))^{-1} \real(Z_i^*\nu)$\;
\State Compute $Z_{i+1} \gets (Z_i+ \eta)/|Z_i + \eta |$\;\Comment{Terminal condition is verified after this step}
\EndFor
\end{algorithmic}
\caption{Schur form Rayleigh quotient iteration for complex tensor eigenpairs}
\label{alg:ComplexEig}
\end{algorithm}
The corresponding Schur form in the real case is different from algorithm $\NCM$ in \cite{JaffeWeissNadler}. The authors found $\ONCM$ outperforms $\NCM$. Performance comparison between our RQI Newton and Schur forms is inconclusive, however, we focus on the Schur form for ease of implementation.

In \cref{alg:ComplexEig}, the initial point $Z_0$ on the unitary sphere $Z_0^*Z_0 = 1$ is chosen randomly. To make sure we only count distinct pairs under the equivalent relation, we keep a table tracking all eigenpairs we have found within our search, identifying pairs where $Z$ is scaled by a $(m-2)$-th root of $1$. This algorithm finds all real pairs if the pairs are isolated. The approach has the advantage of a quadratic convergence algorithm for an individual pair. It does not work well with nonisolated zero, but some remedies could be applied. The speed to find all pairs depends on the distribution of the pairs; generally it works well for randomly generated tensors.

In \href{https://github.com/dnguyend/rayleigh_newton/blob/master/colab/UZPairsEigenTensor.ipynb}{UZPairsEigenTensor.ipynb} in \cite{rayleigh_newton_code} (containing calculations for all examples here), we compare with examples in \cite{CuiDaiNie}. For several examples in that paper, RQI outperforms by a large time factor. For nonisolated eigenpairs or infinite pairs, both approaches need modifications. We will not go into details, but note the case where $\cT$ is the gradient of the Motzkin's polynomial $\hcT(Z^{[6]}) = z_0^4z_1 +z_0^2z_1^4 + z_2^6 - 3z_0^2z_1^2z_2^2$ \cite[example 5.9]{CARTWRIGHT2013942}, we found new complex pairs. The authors found 25 eigenpairs out of the 31 expected pairs. The 6 missing pairs are complex pairs, (satisfying $Z^{\ft}Z=0$, the authors used the E-pair normalization $Z^{\ft}Z=1$) with eigenvalues $\frac{1}{12}$, eigenvectors of the form $(z_0, -\bar{z}_0, 0)$, $4z_0^4=-1$ (2 equivalent pairs) and $\frac{3}{16}$, with eigenvectors $(\pm\frac{\sqrt{-1}}{2}, \pm\frac{\sqrt{-1}}{2}, \frac{\sqrt{2}}{2})$ (4 equivalent pairs). Up to equivalence, we found $6$ real pairs with eigenvalue $0$, while the paper counted 14. A perturbation, adding a small diagonal tensor shows there are indeed 14 ($8$ complex and $6$ real) pairs of small eigenvalues that collapse to the $6$ real pairs, the eigenvectors $(1, 0, 0)$ and $(0, 1, 0)$ each has multiplicity $5$, and each eigenvector $(\pm\frac{1}{\sqrt{3}}, \pm\frac{1}{\sqrt{3}}, \frac{1}{\sqrt{3}})$ has multiplicity $1$ for the count of $14$. See \cite[table 7]{Homotopy} or \cite{rayleigh_newton_code} for the remaining pairs.

\subsubsection{$\cB$ eigenpairs}\label{rem:Btensor}
For the general problem $\bL(\bx, \blbd) = \cT(\bx^{[m-1]}) - \cB(\bx^{[d-1]})\lambda$, it is easy to derive several RQIs based on different constraints, e.g. by using a linear constraint $\bC^{\ft} Z=1$, where $\bC$ is a randomly generated real vector \cite{Homotopy}, thus $\EL=\C$, which we will not go into detail.

We discuss briefly the case where $\cT, \cB$ are symmetric, $\cT = (1/m)\hcT', \cB = (1/d)\hcB'$ where $\hcT, \hcB$ are scalar homogeneous polynomials of order $m$ and $d$, over $\R$. Assuming $\hcB(\bx^{[d]}) = 1$ has a nonempty solution set. We can use the constraint $\hcB(\bx^{[d]}) = 1$ with the rescaling\slash projection retraction discussed in \cref{subsec:rescaleProject}. The tangent space at $\bx$ is defined by the equation $\cB(\bx^{[d-1]})^{\ft}\eta = 0$. We can take  $\bLlbd^-(\bx, \blbd)=\bx^{\ft}$, resulting in the Rayleigh quotient $\cR(\bx) = \hcT(\bx^{[m]})$. The projection is defined by $\Pi(\bx) \omega = \omega - \cB(\bx^{[d-1]})\bx^{\ft}\omega$, hence $\Imag(\Pi(\bx)) =\EE_{\bx}$ is defined by $\bx^{\ft}\omega = 0$. The RQI step can be computed as above.

To analyze cubic convergence, $\cR'(\bx_*)\eta = m \cT(\bx_*^{[m-1]})^{\ft}\eta= m\lambda\cB(\bx_*^{[d-1]})^{\ft}\eta=0$, at an eigenpair $(\bx_*, \cR(\bx_*)$ and $\eta\in T_{\bx_*}$. By \cref{eq:scaleRetract}, $\bR_{\eta\eta}(\bx, 0)\eta^{[2]}$ is proportional to $\bx$,
  $$\Pi(\bx_*)\bLx(\bx_*, \lambda_*)\bx_* =
\Pi(\bx_*)((m-1)\cT(\bx_*^{[m-1]}) - (d-1)\lambda\cB(\bx_*^{[d-1]})) = (m-d)\lambda\Pi(\bx_*)\cB(\bx_*^{[d-1]}) = 0.$$
From \cref{eq:BGRQI}, only $\Pi\bL_{\bx\bx}$ remains in $\bG_{\bL}$. Thus, for $m=d=2$, we have a {\it cubically convergence} iteration for the {\it generalized eigenvalue problem} $A\bx - B\bx\blbd=0$ for symmetric matrices $A$ and $B$ ($\cT(X) = AX, \cB=BX$). 
$$X_{i+1} = \{(A- B\bx_i^{\ft}A\bx_i)^{-1}B\bx_i\}_{rsc}=\{(B^{-1}A- (\bx_i^{\ft}A\bx_i)I)^{-1}\bx_i\}_{rscl}$$
if $B$ is invertible, where $X_{rscl}$ is the rescaling of $X\in\R^n$ to satisfy $\bx^{\ft}B\bx=1$.  Table \ref{tab:num_sim} compares RQI and Rayleigh-Chebyshev iteration for $m=3,d=2$ (see also \cite{rayleigh_newton_code}).
\begin{table}
  \begin{tabular}{c c c c c}
\hline
Iteration&Tensor RQI&Tensor RC &EigC RQI&EigC RC\\
\hline
1&2.594e-04&2.042e-05 &6.986e-06&6.986e-06\\
2&3.080e-07&8.671e-14 &2.368e-11&6.511e-15\\
3&2.509e-13&2.690e-39 &3.500e-22&5.208e-42\\
4&1.049e-25&5.503e-115 &9.982e-44&\\
\hline
\end{tabular}
  \caption{Examples of residual errors for Rayleigh (RQI) and Rayleigh-Chebyshev (RC) iterations for tensor eigenpairs (Tensor, \cref{rem:Btensor}, $n=6, m=3, d=2$) and eigenvalue with a constant term (EigC) (\cref{subsec:eic}, $n=10$).
  }
\label{tab:num_sim}
\end{table}
\subsection{Constrained optimization}\label{subsec:opt}
Recall in this problem, $\bL(\bx, \blbd) = \nabla f(\bx) - \bC'(\bx)^{\ft}\lambda$ for a real-valued function $f$, $\bLlbd(\bx, \blbd) = -\bC'(\bx)^{\ft}$. Let $\bCx = \bC'(\bx)$ then the Moore-Penrose left inverse is $\bLlbd^-(\bx) = -(\bCx\bCx^{\ft})^{-1}\bCx$, and the equation $\bLlbd^-(\bx)\bL(\bx, \blbd) = 0$ implies $\lambda = (\bCx\bCx^{\ft})^{-1}\bCx\nabla f(\bx)$ which is known classically. This implies $\Pi(\bx) = I_{\EE} - \bCx^{\ft}(\bCx\bCx^{\ft})^{-1}\bCx$. Newton's method of this type has been studied before \cite{Goodman,OvertonNocedal}, using the Newton form to solve for $\bx^{\cN_T}$. Let $Q$ be an orthogonal basis of the range of $\Pi(\bx)$, (which is the nullspace of $\bCx$), constructed from the QR decomposition of $\bCx^{\ft}$. Then $Q^{\ft}Q = I$ and $\Pi(\bx) = QQ^{\ft}$. In \cref{eq:NewtonForm}, $Q_{\Pi} = Q_T = Q$, and the Newton increment reduces to
\begin{equation}\bx^{\cN_T} = -Q(Q^{\ft}\bLx(\bx, \blbd)Q)^{-1}Q^{\ft}\nabla f(\bx)\end{equation}
with $Q^{\ft}\nabla f(\bx)$ called the projected gradient, $Q^{\ft}\bLx(\bx, \blbd)Q$ the projected Hessian (the manifold literature calls $\Pi(\bx)\nabla f(\bx)$ the projected gradient and $\Pi(\bx)\bLx(\bx, \blbd)$ the projected Hessian). The Schur complement method is also used, see \cite[Chapter 10]{boyd_vandenberghe_2004}.

\subsection{Eigenvector and eigenvector with a constant term}\label{subsec:eic}
We have shown our RQI with the constraint $\bx^{\ft}\bx -1 =0$ is the classical RQI for the eigenvalue problem, with $\bLx(\bx, \blbd)\eta = (A-\lambda I)\eta$, $\bLlbd(\bx, \blbd) = -\bx$. With $\cR(\bx) = \bx^{\ft} A \bx$, $\cR'(\bx)\eta = \eta^{\ft} A \bx + \bx^{\ft} A\eta = \eta^{\ft}(A + A^{\ft})X$. Thus,
$$\bG_{\bL}(\bx)\eta^{[2]}= 2\bL_{\bx\blbd}(\bx, \lambda; \eta,\cR'(\bx, \eta)) + \bLx(\bx, \blbd)\bR_{\eta\eta}(\bx;\eta^{[2]}) =-2\eta\eta^{\ft}( A + A^{\ft})\bx- (A-\lambda I)\bx|\eta|^2.$$

When $A$ is normal, at an eigenvector $(\bx_*, \lambda)$ we have $A\bx_* = A^{\ft}\bx_*= \lambda \bx_*$, hence $\cR'(\bx_*)\eta = 2\lambda \eta^{\ft}\bx_* =0$, thus $\bG(\bx_*) = 0$, therefore we have cubic convergence.

When $A$ is not normal, at an eigenvector $\bx_*$, only $2\eta^{\ft}A^{\ft}\bx_*$ may be nonzero. To the extent we have tested, the Rayleigh-Chebyshev iteration seems not competitive for generic matrices, as although it  eventually has cubic convergence, its area of convergence seems smaller than RQI.

For nonnormal matrices, the two-sided Rayleigh quotient algorithm in \cite{Ostrowski1959} has cubic convergence. It is a special case of the two-sided nonlinear eigenvalue algorithm in \cref{subsec:nonlinear_rayleigh}.

For $\tb\in\R^n, \tb \neq 0$, an RQI for $\bL(\bx, \lambda) = A\bx -\lambda \bx -\tb=0$ could be derived similarly. $\bLx, \bLlbd, \bLlbdm$ are as above, the Rayleigh quotient is $\lambda = \bx^{\ft}A\bx - \bx^{\ft}\tb$. In the resulting Schur form iteration, $\bLx^{-1}\bL = \bx -(A-\lambda I)^{-1}\tb$. Define $Z_{unit} = Z/\|Z\|$ for $Z\in\R^n$  then
\begin{gather}\eta = -\bx + (A-\blbd I)^{-1}\tb + (A-\blbd I)^{-1}\bx\frac{1-\bx^{\ft}(A-\lambda I)^{-1}\tb}{\bx^{\ft}(A-\lambda I)^{-1}\bx}\label{eq:etaNonLinear}\\
  \bx_{i+1} =((A-\blbd I)^{-1}\tb + (A-\blbd I)^{-1}\bx\frac{1-\bx^{\ft}(A-\lambda I)^{-1}\tb}{\bx^{\ft}(A-\lambda I)^{-1}\bx})_{unit}
\end{gather}
The Chebyshev iteration is more stable than the eigenvector case as $(A-\blbd I)$ is often invertible near a solution. While requiring less iterations, it still does not outperform the RQI in total time. This does not rule out applications to special matrices. The expression for $\bG_{\bL}$ is unchanged, with $\eta$ in \cref{eq:etaNonLinear} and $\cR'(\bx, \eta) =  \eta^{\ft}(A+A^{\ft})\bx - \eta^{\ft}\tb$. The RQI does not converge cubically even for symmetric matrices as $\eta^{\ft} A\bx=\eta^{T}\tb$ near a solution. The Chebyshev increment is
$$\eta^{\bC} = \eta -\frac{1}{2}\{(A-\lambda I)^{-1}\bG_{\bL} - \frac{\bx^{\ft}(A-\lambda I)^{-1}\bG_{\bL}}{\bx^{\ft}(A-\lambda I)^{-1}\bx}(A-\lambda I)^{-1}\bx
\}$$
and the Chebyshev iteration is $\bx_{i+1} = (\bx_i + \eta^{\bC})_{unit}$. The last two columns of \cref{tab:num_sim} shows the RQI and Rayleigh-Chebyshev iterations for this case.
\subsection{Nonlinear eigenvalue problem}
\label{subsec:nonlinear_rayleigh}
A nonlinear RQI appeared in algorithm 4.9 in \cite{guttel_tisseur_2017}, proposed in \cite{schreiber_thesis}. We consider $\EE=\R^n$, $\EL=\R$. Recall $\bL(\bx, \blbd) = \bP(\blbd)\bx$ where $\bP$ is a square matrix with polynomial entries.

Assume the unit sphere constraint $\bC(\bx) = \frac{1}{2}(\bx^{\ft}\bx - 1)$. For $\eta\in \EE$, we have $\bC'(\bx; \eta) = \bx^{\ft}\eta$ and $\bLx(\bx, \lambda)\eta = \bP(\lambda)\eta$ for $\eta\in \EE$, $\bLlbd(\bx, \lambda) \delta = \bP'(\lambda)\bx \delta$ for $\delta\in \R$. Assuming $\bx^{\ft}\bP'(\lambda)\bx\neq 0$, $\bLlbdm(\bx, \lambda)\omega = \frac{\bx^{\ft}}{\bx^{\ft}\bP'(\lambda)\bx}$ is a left inverse of $\bLlbd(\bx, \lambda)$. The equation $\bLlbdm(\bx, \lambda)\bL(\bx, \lambda)=0$ for the Rayleigh quotient implies $\bx^{\ft}\bP(\lambda)\bx = 0$, and $\lambda=\cR(\bx)$ is solved from here. The projection is
$$\Pi(\bx)\omega:= \omega - \bLlbd(\bx, \lambda)\bLlbdm(\bx, \lambda)\omega= \omega -
\bP'(\lambda)\bx\frac{\bx^{\ft}\omega}{\bx^{\ft}\bP'(\lambda)\bx}.$$
Applying \cref{alg:alRQI} in Schur form, $\bLx(\bx, \lambda)^{-1}\bL(\bx, \lambda) = \bx, \bLx(\bx, \lambda)^{-1}\bLlbd(\bx, \lambda)=P(\lambda)^{-1}P'(\lambda)\bx$,
$$\begin{gathered} \eta = \bx^{\cN_T} = -X + P(\lambda)^{-1}P'(\lambda)\bx\frac{\bx^{\ft}\bx}{\bx^{\ft}P(\lambda)^{-1}P'(\lambda)\bx}.
\end{gathered}$$
Thus, $\eta +X$ is proportional to $P(\lambda)^{-1}P'(\lambda)\bx$, and together with the retraction, the iteration could be given as $\bx_{i+1} = \frac{P(\lambda_i)^{-1}P'(\lambda_i)\bx_i}{|P(\lambda_i)^{-1}P'(\lambda_i)\bx_i|}$.
For cubic convergence analysis, from \cref{eq:BGRQI}
\begin{equation}\bG_{\bL}(\bx)\eta^{[2]} = P^{(2)}(\lambda)\bx(\cR'(\bx, \eta))^2 + 2P'(\lambda)\eta\cR'(\bx, \eta) - P(\lambda)\bx|\eta|^2.\label{eq:bGnonlinear}\end{equation}
At a solution, $P(\lambda_*)\bx_*=0$, thus the last term vanishes. From the implicit function theorem
$$\cR'(\bx;\eta) = -\frac{1}{\bx^{\ft} P'(\lambda)\bx}(\eta^{\ft}P(\lambda)\bx + \bx^{\ft} P(\lambda)\eta).$$
If $\bP$ is normal, $\cR'(\bx_*, \eta) = 0$ at a solution, hence $\bG_{\bL}(\bx_*)=0$, implying cubic convergence.

{\it The two-sided iteration}. Consider $\EE = \R^{2n}$ and $\EL = \R^2$, $\bx = \begin{bmatrix}u^{\ft} & v^{\ft}\end{bmatrix}^{\ft}\in \EE$, $u, v\in\R^n$ and $\blambda=[\lambda_1, \lambda_2]^{\ft}\in\EL$. Set $\bC(\bx) = \frac{1}{2}[u^{\ft}u -1, v^{\ft}v-1]^{\ft}\in \R^2$ and
$$\begin{gathered}\hP(\blambda)= \begin{bmatrix}0 & \bP(\lambda_2)^{\ft}\\ \bP(\lambda_1) & 0\end{bmatrix},\quad
\hP(\blambda)^{-1}= \begin{bmatrix}0 & \bP(\lambda_1)^{-1}\\ (\bP(\lambda_2)^{\ft})^{-1} & 0\end{bmatrix},
    \\
  \bL(\bx, \blambda) =\hP(\blambda)\bx =\begin{bmatrix}\bP(\lambda_2)^{\ft}v\\ \bP(\lambda_1)u\end{bmatrix}.
    \end{gathered}$$
In matrix form, $\bL_{\bx}(\bx, \blambda)=\hP(\blambda)$, $\bLlbd$ and a left inverse $\bLlbdm$ are 
\begin{equation}\bL_{\blambda}(\bx, \blambda) = \begin{bmatrix}0 & \bP'(\lambda_2)^{\ft}v\\ \bP'(\lambda_1)u & 0\end{bmatrix},\quad
  \bLlbdm(\bx, \blambda) =  \begin{bmatrix}0 & \frac{v^{\ft}}{v^{\ft}\bP'(\lambda_1)u}\\
    \frac{u^{\ft}}{u^{\ft}\bP'(\lambda_2)^{\ft}v} & 0  \end{bmatrix}\label{eq:nonlinLlbdLlbdm}\end{equation}
  assuming $v^{\ft}\bP'(\lambda_1)u, u^{\ft}\bP'(\lambda_2)v\neq 0$. Then, $\blambda=\cR(\bx)$ satisfies $\bLlbdm(\bx, \blambda)\bL(\bx, \blambda) = 0$
\begin{equation}v^{\ft}\bP(\lambda_1)u = 0,\quad u^{\ft}\bP(\lambda_2)^{\ft}v = 0. \label{eq:Raytwosides} \end{equation}

We have $\bC'(\bx) = \bdiag(u^{\ft}, v^{\ft})\in \R^{2\times 2n}$ where $\bdiag$ denotes a rectangular block diagonal matrix and $\bLx^{-1}\bL =\bx$. In \cref{alg:alRQI}
$$\begin{gathered}
  \bLx^{-1}\bLlbd = \bdiag(\bP(\lambda_1)^{-1}\bP'(\lambda_1)u, (\bP(\lambda_2)^{\ft})^{-1}\bP'(\lambda_2)^{\ft}v)\in \R^{2n\times 2},\\
    \bC'(\bx)\bLx^{-1}\bLlbd = \diag(u^{\ft}\bP(\lambda_1)^{-1}\bP'(\lambda_1)u, v^{\ft}(\bP(\lambda_2)^{\ft})^{-1}\bP'(\lambda_2)^{\ft}v)\in \R^{2\times 2},\\
    \bC'\bLx^{-1}\bL = (u^{\ft}u, v^{\ft}v)\in \R^2,\\
  \eta  = - \bx+ \begin{bmatrix}u^{\ft}u/(u^{\ft}\bP(\lambda_1)^{-1}\bP'(\lambda_1)u)\bP(\lambda_1)^{-1}\bP'(\lambda_1)u,\\
v^{\ft}v/(v^{\ft}\bP(\lambda_1)^{-1}\bP'(\lambda_1)v)    (\bP(\lambda_2)^{\ft})^{-1}\bP(\lambda_2)^{\ft}v\end{bmatrix}
\end{gathered}
$$
where $\eta$ is the Newton increment. Thus, components  of $X+\eta$ are proportional to $\bP(\lambda_1)^{-1}\bP(\lambda_1)u$ and $\bP(\lambda_2)^{\ft})^{-1}\bP(\lambda_2)^{\ft}v$. Hence, the Newton step is
$$u_{i+1} = \frac{(\bP(\lambda_1))^{-1}\bP(\lambda_1)u_i}{|(\bP(\lambda_1))^{-1}\bP(\lambda_1)u_i|}, v_{i+1} = \frac{(\bP(\lambda_2)^{\ft})^{-1}\bP(\lambda_2)^{\ft}v_i}{|(\bP(\lambda_2)^{\ft})^{-1}\bP(\lambda_2)^{\ft}v_i|}.$$
For $\lambda \in \R$, $v^{\ft}\bP(\lambda)u = u^{\ft}\bP(\lambda)^{\ft}v$, so in \cref{eq:Raytwosides} we can choose $\lambda_1=\lambda_2$, the resulting iteration converges to left and right eigenvectors of the same eigenvalue. We have cubic convergence in this case (if $\lambda_1\neq \lambda_2$, the iteration may still converge quadratically to different eigenvalues). In fact, if $\lambda(u, v)$ solves $v^{\ft}\bP(\lambda)u=0$, differentiate in the tangent direction $\eta = (\eta_u, \eta_v)$
$$\begin{gathered}\eta_v^{\ft} \bP(\lambda) u +  v^{\ft} \bP(\lambda) \eta_u +  v^{\ft} \bP'(\lambda) u\lambda'(\bx, \eta) =0\end{gathered}.$$
At an eigenvalue $\lambda_*$ with right and left eigenvectors $u_*, v_*$, the first two terms vanish, hence $\lambda'(\bx_*, \eta)=0$ if $ v^{\ft}_* \bP'(\lambda_*) u_*\neq 0$. For $\bx_* = (u_*, v_*), \blambda_* = (\lambda_*, \lambda_*)$, we have $\bG_{\bL}(\bx_*, \lambda_*) = 0$ since
  $$\bLx(\bx_*, \blambda_*)\bR_{\eta\eta}(\bx_*, 0;\eta^{[2]})=[-(\bP^{\ft}(\lambda_*)v_*|\eta_v|^2)^{\ft},
    -(\bP(\lambda_*)u_*)^{\ft}|\eta_u|^2]^{\ft} = 0.$$
\section{Concluding remarks} We gave an effective procedure to derive Rayleigh quotient iterations for constrained equations with nonlinear multipliers and provided a clear analysis of second and third-order convergence. The theory developed here explains classical cubic convergence results. It also gives an effective algorithm to find all complex tensor eigenpairs. In future work, we will provide a Kantorovich's type theorem \cite{FERREIRA2002}. The non-isolated zeros case could also be studied.
\begin{appendices}
  \section{Various proofs}
  \subsection{Proof of \cref{lem:cl2}}\label{appx:SchurEq}
  By direct substitution, using $P H=0$
$$\begin{gathered}
    P B\eta =  P B(-B^{-1} F +B^{-1}H( DB^{-1}H)^{-1} DB^{-1} F) = -P F + P H( DB^{-1}H)^{-1} DB^{-1} F)
    =-P F,\\
    D\eta =  D(-B^{-1}F +B^{-1}H(DB^{-1}H)^{-1}DB^{-1} F) =  D(-B^{-1} F) + DB^{-1} F=0.\end{gathered}$$
\subsection{A lemma for projection functions}
  \begin{lemma}Let $T_0$ and $\EE$ be two vector spaces, $\Omega \subset T_0$ be an open subset, $\Phi:\Omega\to\Lin(\EE, \EE)$ be a smooth projection function. If $Z\in\Imag(\Phi(\xi))\subset \EE$, hence $\Phi(\xi) Z = Z$, then for $\eta\in T_0$
  \begin{equation}\label{eq:secondFForm}
    \Phi(\xi)\Phi'(\xi; \eta) Z = 0.\end{equation}
\end{lemma}
\begin{proof}Differentiate $\Phi(\xi)^2Z = Z$ in direction $\eta$, we get $\Phi(\xi)\Phi'(\xi; \eta)Z +
  \Phi'(\xi; \eta)\Phi(\xi)Z = 0$. Apply $\Phi(\xi)$ to both sides we get $2\Phi(\xi)\Phi'(\xi; \eta)Z = 0.$
\end{proof}  
\subsection{Proof of \cref{prop:prelimOverdeter}}\label{appx:proofOverdeter}
 The first two equalities in 1) are immediate. Identify $\EE$ and $T_0$ with $\R^n$ and $\R^k$ so we can define transpose, if $\cF^{\natural}$ is defined, then $\Phi(\xi)\cF'(\xi)$ is of rank $k$, and $B(\xi):=(\Phi(\xi)\cF'(\xi))^{\ft}\Phi(\xi)\cF'(\xi)\in \Lin(T_0, T_0)$ is invertible. The equation for $\cF^{\natural}$ could be written $B(\xi)\cF^{\natural}(\xi)\omega = (\Phi(\xi)\cF'(\xi))^{\ft}\Phi(\xi)\omega$ for $\omega\in \EE$. Since $B$ is differentiable, $\xi\mapsto B(\xi)^{-1}$ is differentiable, thus $\cF^{\natural}(\xi)$ is differentiable.

  For \cref{eq:JFnat}, 
  differentiated $\Phi(\xi)\cF'(\xi)\cF^{\natural}(\xi)\omega = \Phi(\xi)\omega$ in $\xi$ in direction $\eta$ then apply $\cF^{\natural}(\xi)$ 
  $$\begin{gathered}\Phi(\xi)\cF'(\xi)(\cF^{\natural})'(\xi; \eta) \omega +
  \Phi'(\xi; \eta)\cF'(\xi) \cF^{\natural}(\xi)\omega
  + \Phi(\xi)\cF^{(2)}(\xi; \eta, \cF^{\natural}(\xi)\omega)
  = \Phi'(\xi;\eta)\omega;\\
(\cF^{\natural})'(\xi; \eta) \omega +
\cF^{\natural}(\xi)  \Phi'(\xi; \eta)\cF'(\xi) \cF^{\natural}(\xi)\omega
  + \cF^{\natural}(\xi)\cF^{(2)}(\xi; \eta, \cF^{\natural}(\xi)\omega)
  = \cF^{\natural}(\xi)\Phi'(\xi;\eta)\omega.
\end{gathered}  $$
  For 2), since $\cF^{\natural}(\xi)\Phi(\xi) = \cF^{\natural}(\xi)$, $\Phi(\xi)Y - \cF(\xi)=0$ implies
  $$\cF^{\natural}(\xi)(Y - \cF(\xi)) = \cF^{\natural}(\xi)[\Phi(\xi)Y - \cF(\xi)] = 0.$$  
  The converse follows by simplifying $\Phi(\xi)\cF'(\xi)\cF^{\natural}(\xi)[Y - \cF(\xi)]=0$, using $\Phi(\xi)\cF(\xi) = \cF(\xi)$.
  
  Consider $G:(\xi, Y) \mapsto \cF^{\natural}(\xi)[Y - \cF(\xi)]$, mapping $D_1\times T_0$ to $T_0$. If $(\xi, Y)$ satisfies $\cF(\xi) = \Phi(\xi)Y$
  $$\begin{gathered}G(., Y)'(\xi; \eta) = (\cF^{\natural})'(\xi; \eta)[Y - \cF(\xi) ] -
    \cF^{\natural}(\xi)\cF'(\xi)\eta
    = (\cF^{\natural})'(\xi; \eta)[Y - \Phi(\xi)Y] -\eta\\ =
    -\eta + \cF^{\natural}(\xi)\Phi'(\xi;\eta)[Y- \Phi(\xi)Y] = -(\eta - \cF^{\natural}(\xi)\Phi'(\xi; \eta)Y) = -(I_{\EE} - \cF^{\natural}(\xi)\Phi'(\xi; ., Y))\eta
  \end{gathered}$$
  for $\eta\in T_0$, where using \cref{eq:JFnat} to expand $(\cF^{\natural})'$, we note $\cF^{\natural}(\xi)[Y - \Phi(\xi)Y]=0$, hence the first and last terms in \cref{eq:JFnat} vanish, then use \cref{eq:secondFForm}. To apply the implicit function theorem for $G$, note that $\cF^{\natural}(\xi)\Phi'(\xi; .; Z)=0$ for $Z= \cF(\xi)$, hence, $I_{T_0} - \cF^{\natural}(\xi)\Phi'(\xi; .; Z) = I_{T_0}$ is invertible, and near $Z$ it is still invertible, thus we have \cref{eq:bJFdiam} from
  $$\begin{gathered}
  \cF_{\diamond}'(Y)\omega = -(G(., Y)')^{-1}(\xi) G(\xi, .)'(Y; \omega) = (I_{T_0} - \cF^{\natural}(\xi)\Phi(\xi; .; Y)')^{-1}\cF^{\natural}(\xi)\omega,\\
 \cF^{\natural}(\xi)\omega = \cF_{\diamond}'(Y)\omega - \cF^{\natural}(\xi)\Phi'(\xi; \cF_{\diamond}'(Y)\omega)Y.
  \end{gathered}$$  
  Differentiate the above in $Y$ in direction $\omega_1\in \EE$, noting $\xi = \cF_{\diamond}(Y)$ (and not constant) we have
  $$\begin{gathered}
  (\cF^{\natural})'(\xi; \cF_{\diamond}'(Y, \omega_1))\omega    
=    
  \cF_{\diamond}^{(2)}(Y; [\omega_1,  \omega]) - \cF^{\natural}(\xi)\Phi'(\xi; \cF_{\diamond}^{(2)}(Y; [\omega_1,  \omega])) Y - \cF^{\natural}(\xi)\Phi'(\xi; \cF_{\diamond}'(Y; \omega)) \omega_1
 \\
-  (\cF^{\natural})'(\xi; \cF_{\diamond}'(Y, \omega_1))\Phi'(\xi; \cF_{\diamond}'(Y; \omega))Y
- \cF^{\natural}(\xi)\Phi^{(2)}(\xi; [\cF_{\diamond}'(Y; \omega_1), \cF_{\diamond}'(Y; \omega)]) Y.
     \end{gathered}$$
  For $Y = \cF(\xi) =: Z, \omega_1= \omega$, using \cref{eq:bJFdiamZ,eq:secondFForm} to simplify the above
  $$\begin{gathered}    
    \cF_{\diamond}^{(2)}(Z; \omega^{[2]}) =    
  (\cF^{\natural})'(\xi; \cF^{\natural}(\xi, \omega))\omega
 + \cF^{\natural}(\xi)\Phi'(\xi; \cF_{\diamond}^{(2)}(Z; \omega^{[2]})) Z + \cF^{\natural}(\xi)\Phi'(\xi; \cF^{\natural}(\xi)\omega) \omega
 \\
+  (\cF^{\natural})'(\xi; \cF^{\natural}(\xi)\omega)\Phi'(\xi; \cF^{\natural}(\xi)\omega))Z
+ \cF^{\natural}(\xi)\Phi^{(2)}(\xi; (\cF^{\natural}(\xi)\omega)^{[2]}) Z\\
= (\cF^{\natural})'(\xi; \cF^{\natural}(\xi, \omega))(\omega+\Phi'(\xi; \cF^{\natural}(\xi)\omega)Z)
 + \cF^{\natural}(\xi)\Phi'(\xi; \cF^{\natural}(\xi)\omega) \omega
+ \cF^{\natural}(\xi)\Phi^{(2)}(\xi; (\cF^{\natural}(\xi)\omega)^{[2]}) Z.
     \end{gathered}
  $$
  Let $\nu = \Phi'(\xi; \cF^{\natural}(\xi)\omega)\cF(\xi) - \cF'(\xi)\cF^{\natural}(\xi)\omega$. Differentiate $\Phi(\xi)\cF(\xi) = \cF(\xi)$ in direction $\cF^{\natural}(\xi)\omega$
  $$\Phi'(\xi;\cF^{\natural}(\xi)\omega )\cF(\xi)+
\Phi(\xi)  \cF(\xi)'\cF^{\natural}(\xi)\omega
  = \cF'(\xi)\cF^{\natural}(\xi)\omega.$$
  Thus, $\nu = -\Phi(\xi)  \cF(\xi)'\cF^{\natural}(\xi)\omega$. By \cref{eq:secondFForm},
$\cF^{\natural}(\xi)\Phi'(\xi; \cF^{\natural}(\xi)\omega)\nu=0$. Let $\psi = (\omega+\Phi'(\xi; \cF^{\natural}(\xi)\omega)Z)$. Use \cref{eq:JFnat} and $\cF^{\natural}(\xi)\psi = \cF^{\natural}(\xi)\omega$ to expand $\cF_{\diamond}^{(2)}(Z; \omega^{[2]})$ further to
$$\begin{gathered}
\cF_{\diamond}^{(2)}(Z; \omega^{[2]})=  
  \cF^{\natural}(\xi)\{-\cF^{(2)}(\xi; (\cF^{\natural}(\xi, \omega))^{[2]}) +
  \Phi'(\xi;\cF^{\natural}(\xi, \omega))[\omega+ \nu] \} +\\
 \cF^{\natural}(\xi)\Phi'(\xi; \cF^{\natural}(\xi)\omega) \omega
 + \cF^{\natural}(\xi)\Phi^{(2)}(\xi; (\cF^{\natural}(\xi)\omega)^{[2]}) Z  \\
 =  \cF^{\natural}(\xi)\{-\cF^{(2)}(\xi; (\cF^{\natural}(\xi)\omega)^{[2]}) +
2\Phi'(\xi; \cF^{\natural}(\xi)\omega) \omega
 + \Phi^{(2)}(\xi; (\cF^{\natural}(\xi)\omega)^{[2]}) Z \}.
\end{gathered}$$
This proves \cref{eq:jFDia2}. Higher derivatives of $\cF_{\diamond}$ could be computed by repeated differentiation of $\cF_{\diamond}'$. Inductively, its derivatives could be expressed algebraically in terms of derivatives of $\cF, \Phi$ and powers of $(I_{T_0} - \cF^{\natural}(\xi)\Phi'(\xi; .)Y)^{-1}$, giving us the $C^k$ condition with bounded derivatives.
\subsection{Proof of \cref{thm:mainOver}}\label{appx:prf:mainOver}
From the mean value theorem, if $\cF'(\xi)$ is bounded by $b_{\cF'}$ in a domain then
$$|\cF(\xi) | = |\cF(\xi) -\cF(\xi_*)|\leq \int_0^1|\cF'(\xi_* + t(\xi - \xi_*))[\xi -\xi_*]| dt \leq b_{\cF'}|\xi - \xi_*|.$$
  For 1), from the Taylor formula of $\cF_{\diamond}$ at $\xi_i$ in 2) of \cref{prop:prelimOverdeter}, set $R_2$ to be the Taylor remainder
  $$\begin{gathered}R_2 =R_2(\xi_*, \xi_i)=\cF_{\diamond}(0) - \cF_{\diamond}(\cF(\xi_i)) - \cF_{\diamond}'(\xi_i)[0-\cF(\xi_i)]= \xi_* - \xi_i - \delta_{2, i}\\
|R_2| \leq \int_0^1
(1-t)|\cF_{\diamond}^{(2)}((1-t)\cF(\xi_i) + t\cF(\xi_*);[\cF(\xi_*)-\cF(\xi_i)]^{[2]})|dt\leq \frac{c_{\cF_{\diamond}}}{2}b_{\cF'}^2|\xi_*-\xi_i|^2\end{gathered}$$
in a ball in $D_1$ where $\cF_{\diamond}$ is defined with bounded derivatives. Boundedness and bilinearity of $\cF_{\diamond}^{(2)}$ implies there is $c_{\cF_{\diamond}}>0$ such that $|\cF_{\diamond}^{(2)}(\xi;\delta^{[2]})| \leq c_{\cF_{\diamond}}|\delta|^2$ for $\delta\in T_0$. Since $\xi_{i+1} = \fS(\xi_i, \delta_{2, i})$,
  \begin{equation}\label{eq:remXiI1}
    \xi_* - \xi_{i+1} = \xi_* - \fS(\xi_i, \delta_{2, i}) = 
(\xi_* - \xi_i -\delta_{2, i}) - \fS_2(\xi_i, \delta_{2, i})= R_2 - \fS_2(\xi_i, \delta_{2, i})\end{equation}
  provided $\delta_{2, i}$ is defined and $|\delta_{2, i}| < \rho_2$. Thus, in $D_1$, $|\delta_2| = |\cF^{\natural}(\xi)\cF(\xi)|  \leq \|\cF^{\natural}(\xi)\||\cF(\xi)|\leq b_{\cF^{\natural}}b_{\cF'}|\xi-\xi_*|$ if $\|\cF^{\natural}(\xi)\|$ and $|\cF'(\xi)|$ are bounded by $b_{\cF^{\natural}}$ and $b_{\cF'}$, respectively. Together, we have
$$|\xi_{i+1}-\xi_*| \leq (\frac{c_{\cF_{\diamond}}b_{\cF'}^2}{2} + c_{\fS,2}(b_{\cF^{\natural}}b_{\cF'})^2)|\xi_i-\xi_*|^2.$$
For $\xi_i$ to be well-defined for all $i$, we also need $|\delta_{2,i}| < \rho_2$. Thus, if $|\xi_0-\xi_*|< \rho_3$ with
$$\rho_3 \leq \min\{\frac{1}{(c_{\cF_{\diamond}}b_{\cF'}^2/2 + c_{\fS,2}(b_{\cF^{\natural}}b_{\cF'})^2)}, \rho_1, \frac{\rho_2}{b_{\cF^{\natural}}b_{\cF'}}\},$$
then $\delta_{2,i}$ and $\xi_i$ are well-defined, $|\xi_i-\xi_*|$ decreases with $i$, $\{\xi_i\}$ converges quadratically to $\xi_*$.

For 2), write $\fS_2(\xi_i, \delta_{2, i}) = \frac{1}{2}\fS(\xi_i, .)^{(2)}(0, \delta_{2, i}^{[2]}) + \fS_3(\xi_i, \delta_{2, i})$ and add $\frac{1}{2}\cF^{\natural}(\xi_*)\cG(\xi_*, \delta_{2,i}^{[2]})=0$ in \cref{eq:remXiI1}
$$\begin{gathered}|R_2 - \fS_2(\xi_i, \delta_{2, i})|
  = |R_2 -\frac{1}{2}\fS(\xi_i, .)^{(2)}(0; \delta_{2, i}^{[2]}) - \fS_3(\xi_i, \delta_{2, i})
  + \frac{1}{2}(\cF^{\natural}(\xi_*)\cF^{(2)}(\xi_*; \delta_{2, i}^{[2]}) + \fS(\xi_*, .)^{(2)}(0; \delta_{2,i}^{[2]} ))|\\
\leq   \int_0^1
(1-t)|\cF_{\diamond}^{(2)}((1-t)\cF(\xi_i) + t\cF(\xi_*); \cF(\xi_i)^{[2]}) -
\cF_{\diamond}^{(2)}(\cF(\xi_*); \cF(\xi_i)^{[2]})
|dt +\\
\frac{1}{2}|\cF_{\diamond}^{(2)}(\cF(\xi_*); \cF(\xi_i)^{[2]}) +\cF^{\natural}(\xi_*)\cF^{(2)}(\xi_*; \delta_{2, i}^{[2]})| + \frac{1}{2}|\fS(\xi_*, .)^{(2)}(0; \delta_{2,i}^{[2]})
-\fS(\xi_i, .)^{(2)}(0; \delta_{2, i}^{[2]})|
+|\fS_3(\xi_i, \delta_{2, i})| \\
\leq \frac{c_{\cF_{\diamond},3}}{3}|\cF(\xi_i)||\cF(\xi_i)|^2 +
\frac{Q}{2} +\frac{c_{\fS;2*}}{2}|\xi_*-\xi||\delta_{2,i}^2 +
c_{\fS, 3}|\delta_{2, i}|^3\\
\end{gathered}$$
if we bound $\cF_{\diamond}^{(3)}(., [\delta_1, \delta_2, \delta_3]) \leq c_{\cF_{\diamond, 3}}|\delta_1||\delta_2||\delta_3|$, hence the integral by $c_{\cF_{\diamond, 3}}\int_0^1(1-t)^2|\cF(\xi_i)||\cF(\xi_i)|^2dt$ using the mean value theorem for $\cF_{\diamond}^{(2)}$ and $\cF(\xi_*) = 0$, and similarly for the terms involving $\fS_3$ and $|\fS(\xi_*, .)^{[2]} - \fS(\xi_i, .)^{[2]}|$, eventually get an expression of the form $c|\xi_*-\xi_i|^3$. It remains to estimate the expression $Q$ below using \cref{eq:jFDia2}, adding $-2\cF^{\natural}(\xi_i)\Phi'(\xi_i, \cF^{\natural}(\xi_i)\cF(\xi_i))\cF(\xi_i)=0$ from \cref{eq:secondFForm}
$$\begin{gathered}Q = |\cF_{\diamond}^{(2)}(\cF(\xi_*); \cF(\xi_i)^{[2]}) +\cF^{\natural}(\xi_*)\cF^{(2)}(\xi_*; \delta_{2, i}^{[2]})|]\leq  | \Phi^{(2)}(\xi_*; (\cF^{\natural}(\xi_*)\cF(\xi_i))^{[2]})\cF(\xi_*)\} |\\
    + |\cF^{\natural}(\xi_*)\cF^{(2)}(\xi_*; (\cF(\xi_i)^{\natural}\cF(\xi_i))^{[2]})
    -\cF^{\natural}(\xi_*)\cF^{(2)}(\xi_*; (\cF(\xi_*)^{\natural}\cF(\xi_i))^{[2]})|
      \\
+2|\cF^{\natural}(\xi_*)\Phi'(\xi_*, \cF^{\natural}(\xi_*)\cF(\xi_i))\cF(\xi_i)
  -\cF^{\natural}(\xi_i)\Phi'(\xi_i, \cF^{\natural}(\xi_i)\cF(\xi_i))\cF(\xi_i)|
\end{gathered}$$
The first term is zero since $\cF(\xi_*) =0$, and the other terms are bounded by the mean value theorem, for example, for the last term, consider the mean value of
$f(t) = \cF^{\natural}(\xi_*)\Phi'(\xi_t, \cF^{\natural}(\xi_t)\cF(\xi_i))\cF(\xi_i)$ with $\xi(t) = \xi_i +t(\xi_*-\xi_i)$, which gives us a term $(\xi_*-\xi_i)$, together with the two terms $\cF(\xi_i)$ and boundedness of derivatives giving an $O(|\xi_*-\xi_i|^3)$ estimate.

For 3), with $R_3$ denotes the third-order Taylor remainder of $\cF_{\diamond}$ centered at $\cF(\xi_i)$, from \cref{eq:bJFdiamZ,eq:jFDia2}, (the term with $2\cF^{\natural}\Phi'\cF$ is zero by \cref{eq:secondFForm})
$$\begin{gathered}R_3 = \cF_{\diamond}(0) - \cF_{\diamond}(\cF(\xi_i)) - (\cF_{\diamond})'(\cF(\xi_i); 0-\cF(\xi_i)) - \frac{1}{2}(\cF_{\diamond})^{(2)}(\cF(\xi_i); [0-\cF(\xi_i)]^{[2]})=\\
\xi_* - \xi_i +\cF^{\natural}(\xi_i)\cF(\xi_i) -\frac{1}{2}
\cF^{\natural}(\xi_i)\{-\cF^{(2)}(\xi_i; (\cF^{\natural}(\xi)\cF(\xi_i))^{[2]}) + \Phi^{(2)}(\xi_i; (\cF^{\natural}(\xi_i)\cF(\xi_i))^{[2]})\cF(\xi_i)\}
\end{gathered}$$
Let $\fS_3$ be the third Taylor remainder for $\fS$, expand $\delta_{3, i}$, add and subtract $\frac{1}{2}\cF^{\natural}(\xi_i)\Phi^{(2)}(\xi_i, \delta_{2, i}^{[2]})\cF(\xi_i)$
$$\begin{gathered}\xi_* - \xi_{i+1} = \xi_* - \xi_i -\delta_{3, i} - \frac{1}{2}\fS^{(2)}(\xi_i, .)(0; \delta_{3, i}^{[2]}) - \fS_3(\xi_i, \delta_{3, i})=\\
  \xi_* - \xi_i -\cF^{\natural}(\xi_i)\cF(\xi_i) +\frac{1}{2}\cF^{\natural}(\xi_i)\cF^{(2)}(\xi_i; \delta_{2, i}^{[2]})+\frac{1}{2}\fS(\xi_i, .)^{(2)}(0; \delta_{2, i}^{[2]})
  \\+ \frac{1}{2}\cF^{\natural}(\xi_i)\Phi^{(2)}(\xi_i, \delta_{2, i}^{[2]})\cF(\xi_i) 
   - \frac{1}{2}\cF^{\natural}(\xi_i)\Phi^{(2)}(\xi_i, \delta_{2, i}^{[2]})\cF(\xi_i)
  - \frac{1}{2}\fS(\xi_i, .)^{(2)}(0; \delta_{3, i}^{[2]}) - \fS_3(\xi_i, \delta_{3, i})\\
  = R_3(\xi_*, \xi_i)-  \frac{1}{2}\cF^{\natural}(\xi_i)\Phi^{(2)}(\xi_i, \delta_{2, i}^{[2]})\cF(\xi_i) + \frac{1}{2}(
\fS'(\xi_i, .)(0; \delta_{2, i}^{[2]})
- \fS(\xi_i, .)'(0; \delta_{3, i}^{[2]}))
  - \fS_3(\xi_i, \delta_{3, i}).
\end{gathered}$$
We show each term of the last expression above is $O(|\xi_i - \xi_*|^3)$. For $R_3$, in a sufficiently small ball, there is a constant $c_{\cF_{\diamond}, 3}$ such that $|\cF_{\diamond}^{(3)}(\xi, \delta^{[3]})| \leq c_{\cF_{\diamond}, 3}|\delta|^{3}$ for $\delta\in T_0$ and
$$|R_3(\xi_*, \xi_i)|\leq \int_0^1\frac{(1-t)^2}{2}c_{\cF_{\diamond, 3}}|\cF(\xi_i)|^{3}dt \leq \frac{1}{6}c_{\cF_{\diamond, 3}}b_{\cF'}^3|\xi_i -\xi_*|^{3}.$$
For the next term, bound $\|\cF^{\natural}(\xi_i)\|$ by a constant in a ball, $\Phi^{(2)}(\xi, \delta_{2, i}^{[2]})$ by $b_{\Phi,2}|\delta_{2, i}|^2$ for a constant $b_{\Phi,2}$, and $\cF(\xi_i)$ by $b_{\cF'}|\xi_i-\xi_*|$. Since $\fS(\xi_i, .)^{(2)}$ is bilinear, expand $\delta_{3,i} = \delta_{2, i} + (\delta_{3,i}-\delta_{2, i})$, noting $|\delta_{2, i}| < c_3 |\xi_*-\xi_i|$,
$|\delta_{3,i}-\delta_{2, i}| =\frac{1}{2}|\cF^{\natural}(\xi)\cG(\xi_i; \delta_{2, i}^{[2]}) \leq c_{4*}|\delta_{2, i}|^2\leq c_4|\xi_*-\xi_i|^2$ for constants $c_3, c_4$, the third term is bounded by
$$|\fS(\xi_i, .)^{(2)}(0; [\delta_{2, i},\delta_{3, i}-\delta_{2, i}])|
+\frac{1}{2}|\fS(\xi_i, .)^{(2)}(0; [\delta_{3, i}-\delta_{2, i}]^{[2]}|
\leq c_{\fS, 2}c_{3}c_4|\xi_*-\xi_i|^3 + \frac{c_{\fS, 2}c^2_4}{2}|\xi_*-\xi_i|^4.$$
The $\fS_3$ term is third order. The choice of the initial point for the series to be well-defined is similar to the quadratic case.
\end{appendices}
\bibliographystyle{elsarticle-num}
\bibliography{RayleighNewton}
\end{document}